\documentclass[final, 3p, times]{elsarticle}

\usepackage{lineno, hyperref, float, amsmath}
\usepackage{amsfonts, amssymb, color, comment, graphicx, subcaption}
\usepackage[dvipsnames]{xcolor}

\graphicspath{ {Images/} }
\hypersetup{colorlinks = true}
\AtBeginDocument{\hypersetup{citecolor = blue}}

\journal{CAMWA}

\begin{document}

\begin{frontmatter}

    \title{Efficient adaptive step size control for exponential integrators}

\author{Pranab Jyoti Deka} 
\ead{pranab.deka@uibk.ac.at}

\author{Lukas Einkemmer}
\ead{lukas.einkemmer@uibk.ac.at}

\address{Department of Mathematics, University of Innsbruck, A-6020 Innsbruck, Austria}

%%%%%%%%%%%%%%%%%%%%%%%%%%%%%%%%%%%%%%%%%%%%%%%%%%%%%%%%%%%%%%%%%%%%%%%%%%%%%%%%%%%%%%%%%%%%%%%%%%%%%%%%%%%

\begin{abstract}
Traditional step-size controllers make the tacit assumption that the cost of a time step is independent of the step size. This is equitable with explicit and implicit integrators that use direct solvers. However, in the context of exponential integrators, an iterative approach such as the Krylov method or polynomial interpolation is often employed to compute the action of the required matrix functions. This renders the assumption of constant cost for any given step size invalid. This is a problem for higher-order exponential integrators, as they can take relatively large step sizes based on accuracy considerations. In this manuscript, we consider an adaptive step-size controller for exponential Rosenbrock methods that determines the step size based on the premise of minimizing computational cost. The largest allowed step size, given by accuracy considerations, merely acts as a constraint.  We test this approach on a range of nonlinear partial differential equations. Our results show significant improvements (up to a factor of 4 reduction in the computational cost) over the traditional step-size controller for a wide range of tolerances. 
\end{abstract}


\end{frontmatter}

%%%%%%%%%%%%%%%%%%%%%%%%%%%%%%%%%%%%%%%%%%%%%%%%%%%%%%%%%%%%%%%%%%%%%%%%%%%%%%%%%%%%%%%%%%%%%%%%%%%%%%%%%%%

\section{Introduction}

Solving time dependent partial differential equations (PDEs) numerically is important in almost all fields of science and engineering. Consequently, improvements in numerical algorithms have contributed greatly to better understand a range of natural phenomena and such methods are essential in many industrial settings. Faster numerical methods, in this context, allow us to perform simulations with increased fidelity, e.g.~increasing the number of grid points, including more physical effects, etc. 

While explicit numerical methods are suitable for some problems, for many PDEs a large efficiency improvement can be attained by using implicit time integrators. Consequently, such methods have attained much interest and many software packages have been written to facilitate the use of such methods by practitioners, see e.g.~\cite{hairerII,cvode}. More recently, so-called exponential Rosenbrock integrators have been introduced \cite{Hochbruck09}. We refer to the review article \cite{Ostermann10} for more details. This class of methods linearizes the partial differential equation and then treats a matrix function representing the linear part using Krylov iteration, Leja interpolation, or Taylor methods. Similar to implicit integrators, exponential Rosenbrock methods can take much larger time steps than explicit methods. However, the fact that such methods do not approximate the linear part of the equation (except for the error in the iterative scheme) allows them, in many situations, to take even larger time steps. Moreover, these methods do not suffer from the dichotomy between good behavior on the negative real axis (where the stability function is expected to decay) and on the imaginary axis (where the stability function should have unit magnitude) that afflicts implicit integrators. Because of this, exponential integrators have been used extensively and demonstrated to be superior compared to implicit methods in several situations, see e.g.~\cite{hochbruck1998, Tokman13, Caliari09, einkemmer2017, crouseilles2018, crouseilles2020}. 

To facilitate the use of software packages based on these integrators by practitioners, it is desirable to free the user from explicitly choosing the time step size. Ideally, the user would only prescribe a tolerance and the numerical algorithm would then select an appropriate step size. This can be done using automatic step-size controllers coupled with an error estimator. Ideally, this also frees the user from checking the accuracy of the simulation. For many problems, step-size controllers can also improve computational efficiency by varying the step size as the simulation evolves in time.

Almost all widely used step-size controllers make the assumption that increasing the step size results in a decrease in computational cost. Thus, the step size is chosen such that the error precisely matches the tolerance specified by the user (in practical implementations often a safety factor is imposed to avoid frequent step size rejection). This is a reasonable assumption for explicit Runge\textendash Kutta methods, where the computational cost is independent of the step size. However, this approach is also used in many implicit and exponential software packages. For example, the implicit RADAU5 code \cite[Chap. IV.8]{hairerII} and the implicit multistep based CVODE code \cite{cvode} use this approach. These implicit (or exponential) methods require an iterative solution of a linear system or the iterative computation of the action of certain matrix functions. This is commonly done by iterative methods. However, the number of iterations depends on the spectrum of the matrix. Changing the time step size scales the spectrum and thus also alters the number of iterations. Since the corresponding relationship is not linear, reducing the time step size below what is dictated by accuracy considerations can actually result in an increase in performance, thereby invalidating the assumption that the cost of a time step is independent of the step size.

None of the widely employed step-size controllers are able to exploit this fact. This is problematic for two reasons. First, it reduces the computational efficiency by taking time steps size that do not yield optimal performance. Second, such step-size controllers often do not show a monotonous increase in cost as the tolerance decreases. Thus decreasing the tolerance can actually (sometimes drastically) reduce the run-time required for the simulations. Such behaviour is observed in a range of test problems \cite{Tokman13,luan2017,hochbruck1998} as well as for more realistic physical models \cite{einkemmer2017,blom2016,narayanamurthi2017}. The problem with this behaviour is that the user has to tune the parameters of the method in order to obtain optimal efficiency (i.e.~decreasing the tolerance until the run-time is minimized). Thus, this largely negates the utility of an automatic step-size controller. This behaviour can be observed for exponential integrators as well as implicit Runge--Kutta methods, BDF methods, and implicit-explicit (IMEX) methods.

While all of the considerations made above are valid for implicit schemes just as well as for exponential integrators, the issues raised become even more important for exponential integrators. The reason being that exponential integrators, especially for problems where nonlinear effects are relatively weak, are often able to take even larger time steps than implicit integrators. Thus, exponential integrators when used in conjunction with a traditional step-size controller are more likely to operate in a regime that is problematic.

In the context of ordinary differential equations, the significance of considering a time step size dependent cost function has been recognized in \cite{gustafsson1997}.  In this work, analytically derived cost estimates are employed. However, obtaining a good a priori estimate of the cost is often extremely difficult for (especially nonlinear) PDEs. In \cite{caliari2016} a backward error analysis is used to determine an appropriate step size. However, this approach requires certain information on the spectrum of the matrix, information that is not easily obtained in e.g.~a matrix free implementation, that in nonlinear problems changes each time step. In addition, it is well known that the number of iterations is overestimated and early truncation still happens for many classes of matrices.

In \cite{Einkemmer18}, an adaptive step-size controller has been introduced that explores the space of admissible step sizes (i.e.~step sizes that satisfy the tolerance) dynamically during the simulation and adapts the step size based on the measured cost. This has the advantage that no prior estimates of the cost are needed. In fact, no information of the iterative scheme used or the hardware where the simulation is run, enters the algorithm. Only the computational cost of the previously conducted time steps is used. It was shown in \cite{Einkemmer18} that, for a number of implicit Runge--Kutta methods, this approach reduces the overall computational cost significantly and results in a monotonic relation of the computational effort with the run time.

The goal of the this paper is to consider an adaptive step-size controller for exponential integrators and to investigate its performance. This controller is an extension of the method described in \cite{Einkemmer18} to exponential integrators. As mentioned above, using adaptive step size control is particularly important in the case of exponential integrators. We demonstrate that the developed controller performs well, that is, it increases computational efficiency and removes the non-monotonous behaviour observed in the traditional approach to step size control. We also compare the performance of the adaptive step-size controller to the implicit approach proposed in \cite{Einkemmer18} and find that the present approach can yield improvements in performance of up to an order of magnitude.

The paper is structured as follows. An introduction to exponential integrators and our implementation is presented in section \ref{sec:exp_int}. In section \ref{sec:cost_control}, the principle of the proposed step-size controller is presented. The performance of this step-size controller is then analyzed for some nonlinear problems in section \ref{sec:nl_eqs}. We conclude our study in section \ref{sec:conclude}.

%%%%%%%%%%%%%%%%%%%%%%%%%%%%%%%%%%%%%%%%%%%%%%%%%%%%%%%%%%%%%%%%%%%%%%%%%%%%%%%%%%%%%%%%%%%%%%%%%%%%%%%%%%%

\section{Exponential Integrators}
\label{sec:exp_int}

In this section, we provide an introduction to exponential integrators and Leja interpolation that we use to compute the action of the resulting matrix-vector products. We refer the reader to \cite{Ostermann10} for more details. Let us consider the initial value problem

\begin{equation}
    \frac{\partial u}{\partial t} = f(u), \qquad u(t = 0) = u^0,
    \label{eq:nl_pde}
\end{equation}
where $u \equiv u(x, t)$ in 1D, $u \equiv u(x, y, t)$ in 2D, and $f(u)$ is some nonlinear function of $u$ (usually depends on spatial derivatives of $u$). Linearizing Eq. \ref{eq:nl_pde} about $u^n$, the starting point for a given time step, we get
\begin{equation}
    \frac{\partial u}{\partial t} = \mathcal{J}(u^n) \, u + \mathcal{F}(u), \nonumber
\end{equation}
where $\mathcal{J}(u)$ is the Jacobian of the nonlinear function $f(u)$ and $\mathcal{F}(u) = f(u) - \mathcal{J}(u) \, u$ is the nonlinear remainder. We use exponential Rosenbrock (EXPRB) integrators \citep{Hochbruck06} to solve equations of this form. The simplest of the EXPRB integrators, known as the \textit{exponential Rosenbrock--Euler} integrator, is given by

\begin{equation}
    u^{n + 1} = u^n + \Delta t \varphi_1(\mathcal{J}(u^n) \Delta t) f(u^n),
    \label{eq:exprb2}
\end{equation}
where the superscripts $n$ and $n + 1$ indicate the time steps. The $\varphi_l(z)$ functions are defined by the recursive relation

\begin{equation}
	\varphi_{l + 1}(z) = \frac{1}{z} \left(\varphi_l(z) - \frac{1}{l!} \right), \quad l, \nonumber \geq 1
\end{equation}
with
\begin{equation}
    \varphi_0(z) = e^z, \nonumber
\end{equation}
which corresponds to the matrix exponential. The exponential Rosenbrock--Euler integrator is second-order accurate and only needs the action of one matrix function per time step. An error estimator for Eq. \ref{eq:exprb2} has been developed by \cite{Caliari09}. For the second-order accuracy to hold, it is crucial that the Jacobian is used. In fact, integrators that replace the Jacobian by an arbitrary linear operator, require more stages to obtain a given order. These methods are referred to as either exponential Runge--Kutta integrators or exponential time differencing integrators. We will only consider exponential Rosenbrock integrators in this paper. Many higher order variants of this idea are available in the literature, see e.g.~\cite{Ostermann10, Luan16, Luan18}. Of particular interest in this work are higher order embedded schemes, similar to embedded Runge--Kutta methods, that are a pair of exponential integrators with same internal stages but different order. The difference between these two solutions is then used to cheaply obtain an error estimate for adaptive step-size control. 

In this work, we use the fourth-order (EXPRB43) integrator with a third-order error estimator, presented in \cite{Hochbruck09} and the Butcher tableau of which can be found in \cite{Ostermann10}. The two internal stages are given by $a_n$ and $b_n$, and the third and fourth-order solutions are given by $u_3^{n+1}$ and $u_4^{n+1}$ (Eq. \ref{eq:exprb43}), respectively. The difference between these two solutions gives an error estimate of order three.

\begin{align}
	a^n & = u^n + \frac{1}{2} \varphi_1\left(\frac{1}{2} \mathcal{J}(u^n)  \Delta t \right) f(u^n)  \Delta t \nonumber  \\
	b^n & = u^n + \varphi_1\left(\mathcal{J}(u^n) \Delta t \right) f(u^n) \Delta t \nonumber + \varphi_1\left(\mathcal{J}(u^n)  \Delta t \right) (\mathcal{F}(a^n) - \mathcal{F}(u^n)) \Delta t \nonumber \\
	u_3^{n + 1} & = u^n + \varphi_1\left(\mathcal{J}(u^n) \Delta t\right) f(u^n) \Delta t \nonumber + \varphi_3(\mathcal{J}(u^n) \Delta t)(-14 \mathcal{F}(u^n) + 16 \mathcal{F}(a^n) - 2 \mathcal{F}(b^n)) \Delta t \\
	u_4^{n + 1} & = u_3^{n+1} + \varphi_4(\mathcal{J}(u^n) \Delta t)(36 \mathcal{F}(u^n) - 48 \mathcal{F}(a^n) + 12 \mathcal{F}(b^n)) \Delta t
	\label{eq:exprb43}
\end{align}

%%% ----------------------------------------------------------- %%%
%%% ----------------------------------------------------------- %%%

\subsection{Leja interpolation}

The main computational effort required in an exponential integrator is to evaluate the action of the matrix functions $\varphi_l$. Similar to the treatment of linear solves in implicit schemes, iterative methods are commonly used to treat the large matrices resulting from the spatial discretization of PDEs. Krylov subspace methods, methods based on polynomial interpolation, and Taylor methods are the most common options. A comparison of many of these methods has been conducted in \cite{Bergamaschi06, Caliari14}. In this work, we will exclusively use interpolation at Leja points. However, the developed adaptive step-size controller is expected to work equally well for other strategies.

An effective way of computing the action the matrix exponential and the $\varphi_l$ functions is interpolation at Leja points \cite{Caliari04}. Assuming $\mathbb{K}$ is a compact set and $\mathbb{K} \subset \mathbb{C}$, where $\mathbb{C}$ is the complex plane, a set of Leja points, denoted by $z$ in the following, can be defined recursively as \[\prod_{k = 0}^{j - 1} |z_j - z_k| = \text{max}_{z \in \mathbb{K}}\prod_{k = 0}^{j - 1} |z - z_k|,\] where $z \in \mathbb{K}$ and $j = 1, 2, 3 \hdots$. Conventionally, $|z_0|$ is chosen to be $\text{max}_{z \in \mathbb{K}} |z|$. These recursively defined points were initially studied by \cite{Edrei1940, Leja1957}. Details on the sequence of Leja points and its potential practical applications can be found, for example, in \cite{Reichel90, Baglama98, Caliari04}.

In this work, we approximate the action of the $\varphi_l$ functions by interpolating them as a polynomial on Leja points. The preference for Leja points over the well-known Chebyshev points can be attributed to the fact that the interpolation of a polynomial at $m + 1$ Chebyshev nodes necessitates the re-computation of the $\varphi_{l}$functions/matrix-vector products at the previously computed $m$ nodes. However, Leja points can be generated in a sequence: using $m + 1$ Leja points needs only one extra computation, and the computation at the previous $m$ nodes can be reused.

Here, we present a synopsis of the algorithm that we use in our implementation (following \cite{Caliari14}). Appropriately placing the interpolation points requires the spectral properties of the matrix. Let us suppose that the eigenvalues of the matrix $A$ satisfy
\begin{equation}
    \alpha \leq \mathrm{Re} \; \sigma(A) \leq \nu \leq 0 , \qquad -\beta \leq \mathrm{Im} \; \sigma (A) \leq \beta, \nonumber
\end{equation}
where $\sigma(A)$ denotes the spectrum of $A$; $\alpha$ and $\nu$ are the smallest and largest real eigenvalues respectively, and $\beta$ is the largest, in modulus, imaginary eigenvalue. The values of $\alpha$, $\nu$, and $\beta$ can be obtained by Gershgorin's disk theorem.

One can then construct an ellipse, with semi-major axis $a$ and semi-minor axis $b$, consisting of all the eigenvalues of the matrix $A$. For real eigenvalues, let $c$ be the midpoint of the ellipse and $\gamma$ be one-fourth the distance between the two foci of the ellipse. For the matrix exponential, we interpolate the function $\exp(c + \gamma \xi)$ on pre-computed Leja points ($\xi$) in the interval $[-2, 2]$. 

The $n^\mathrm{th}$ term of the interpolation polynomial $p(z)$ is defined as
\begin{align}
    p_n(z) & = p_{n - 1}(z) + d_n \, y_{n - 1}(z), \nonumber \\
    y_n(z) & = y_{n - 1}(z) \times \left(\frac{z - c}{\gamma} - \xi_n \right), \nonumber
\end{align}
where the $d_i$ correspond to the divided differences of the function $\exp(c + \gamma \xi)$. For imaginary eigenvalues, one can interpolate the function $\exp(c + \gamma \xi)$ on the interval $i[-2, 2]$. To interpolate $\varphi_l$ functions on Leja points, one can simply replace $\exp(c + \gamma \xi)$ with $\varphi_l(c + \gamma \xi)$.

Let us note that the number of Leja points needed for a certain stage, of a multi-stage integrator, to converge depends on a multitude of factors. This includes the step size, the spectrum of the Jacobian, the norm of the function being interpolated, the coefficients of the integrator, the $\varphi_l$ function, and the user-defined defined tolerance. For example, $\varphi_1\left(\frac{1}{2} \mathcal{J}(u^n) \Delta t \right) f(u^n)  \Delta t$ and $\varphi_1\left(\mathcal{J}(u^n) \Delta t \right) f(u^n) \Delta t$ would take different number of Leja points to converge. The coefficients of the polynomial, that are computed using the divided differences algorithm, depend on integrator coefficients (here, $1/2$ and $1$). This results in a nonlinear dependence of these integrator coefficients on the polynomial coefficients, which, in turn, determine, in part, how many Leja points are needed for convergence. Next, the norm of the function, $f(u)$, is usually much larger than that of the nonlinear remainders [$(-14 \mathcal{F}(u^n) + 16 \mathcal{F}(a^n) - 2 \mathcal{F}(b^n))$ and $(36 \mathcal{F}(u^n) - 48 \mathcal{F}(a^n) + 12 \mathcal{F}(b^n))$]. As such, $f(u)$ takes more Leja points to converge than the nonlinear remainders. Another factor that comes into consideration is the order of the $\varphi_l$ function: higher-order $\varphi_l$ functions tend to converge faster, i.e. they need fewer Leja points. Finally, if the step sizes are varied at every time step, which is the case in this study, the number of Leja points used are expected to vary. Even for constant step sizes, the convergence time may vary at every time step owing to the variations in the spectrum of the Jacobian and the function to be interpolated.

%%%%%%%%%%%%%%%%%%%%%%%%%%%%%%%%%%%%%%%%%%%%%%%%%%%%%%%%%%%%%%%%%%%%%%%%%%%%%%%%%%%%%%%%%%%%%%%%%%%%%%%%%%%

\section{Adaptive step-size controller}
\label{sec:cost_control}

\begin{figure*}[t]
    \centering
	\includegraphics[width = 0.6\textwidth]{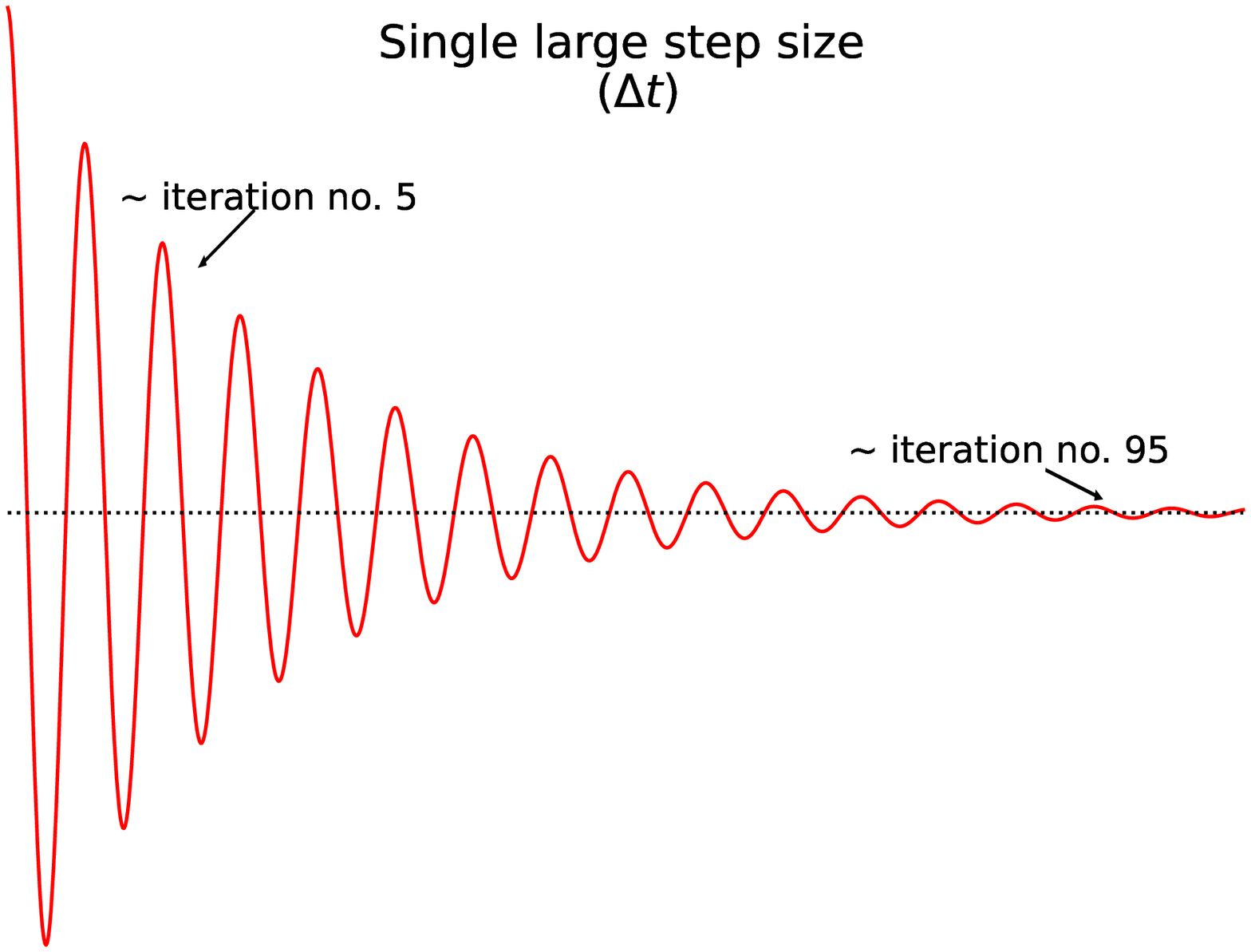} \\
	\includegraphics[width = \textwidth]{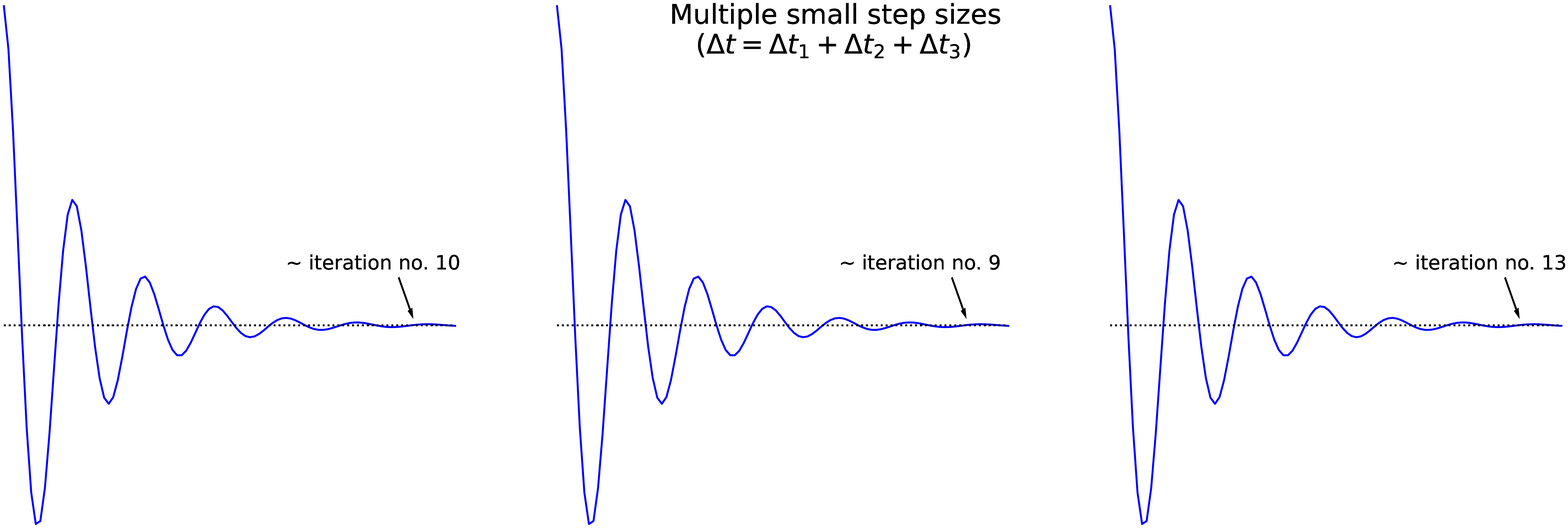}
    \caption{An illustration (not to scale) of convergence as a function of number of iterations for a single large step size versus multiple small step sizes. The number of iterations needed for convergence depends nonlinearly on the step size. This is expected to result in substantially reduced computational cost for multiple small step sizes over one large step size.}
    \label{fig:schematic_dt}
\end{figure*}

To conduct automatic step-size control, an error estimate is essential to ensure that the local error is below the user-specified tolerance. Embedded integrators, that share the internal stages, can be efficiently used as error estimators (with only a small increase in the computation cost). Richardson extrapolation is one of the other commonly used error estimator.

The widely used traditional step-size controller uses the largest possible step size (with a safety factor) that satisfies the prescribed tolerance. This implicitly assumes that the cost of each step size is independent of the step size $\Delta t$. This is true for explicit methods or implicit methods that that solve the corresponding linear systems using direct methods. Let us suppose that the error incurred in the $n^\mathrm{th}$ time step is $e^n = D(\Delta t^n)^{(p + 1)}$, where $p$ is the order of the method used and $D$ is some constant. The tolerance specified by the user is $\mathrm{tol}$.  The optimal step size, for $(n + 1)^\mathrm{th}$ time step, is given by $\mathrm{tol} = D(\Delta t^{n + 1})^{(p + 1)}$. Eliminating $D$, we get

\begin{equation}
    \Delta t^{n+1} = \Delta t^n \times \left(\frac{\mathrm{tol}}{e^n} \right)^{1/(p + 1)}. \nonumber
\end{equation}
Such local step-size controllers are widely used in many different time integration software packages, e.g.~in the RADAU5 code \cite[Chap. IV.8]{hairerII} and the multistep based CVODE code \cite{cvode, eckert2004}. For a mathematical analysis of such methods, we refer the reader to \cite{gustafsson1988,gustafsson1994,soderlind2002,soderlind2006,gustafsson1988}.

For iterative methods (in the context of implicit or exponential time integrators), the computational cost depends on the step size; the larger the step size, the larger the number of iterations needed for the integration to converge (simplistic visual representation in Fig. \ref{fig:schematic_dt}). As such, it is not always beneficial to choose the largest possible step size. Taking this into consideration, \cite{Einkemmer18} developed an adaptive step-size controller where the step size is chosen based on the computational expenses at the previous time steps. This step-size controller is engineered to select step sizes to minimize the computational cost (which might be substantially smaller than the one yielded by the traditional controller). 

This approach works as follows: the step size is adjusted in accordance with the computational cost (c) per unit time step
\begin{equation}
    c^n = \frac{i^n}{\Delta t^n}, \nonumber
\end{equation}
where $i^n$ is the runtime or a proxy, such as the number of matrix-vector products needed in that time step. The goal of this step-size controller is to adjust the step size such that $c^n \longrightarrow \text{min}$. We consider the logarithm of the step size $T^n = \mathrm{ln} \, \Delta t^n$ and the computational cost $C^n(T^n) = \mathrm{ln} \, c^n(\Delta t^n)$. One-dimensional gradient descent is implemented to estimate $T^{n+1}$
\begin{equation}
    T^{n+1} = T^n - \gamma \nabla C^n(T^n), \nonumber
\end{equation}
where $\gamma$ is the learning rate. The gradient can be approximated by taking finite differences
\begin{equation}
    \nabla C^n(T^n) \approx \frac{C^n(T^n) - C^n(T^{n-1})}{T^n - T^{n-1}} \nonumber
\end{equation}
This implies that we are not allowed to choose a constant time step size, i.e.~$T^n \neq T^{n-1}$, as this would not provide any information on how the time step should be changed to optimize the performance. It is worth noting that $C^{n-1}(T^n)$ corresponds to the cost of a step size ($\Delta t^n$) starting from $t^{n-1}$ whereas $C^{n}(T^n)$ is the cost incurred for the same step size ($\Delta t^n$) starting from $t^n$. During the time integration of a problem, we automatically obtain $C^{n-1}(T^{n-1})$ and not $C^{n}(T^{n-1})$. Therefore, we can write the gradient as 

\begin{align}
    \nabla C^n(T^n) & \approx \frac{C^n(T^n) - C^n(T^{n-1})}{T^n - T^{n-1}} \nonumber \\
                    & = \frac{C^n(T^n) - C^{n-1}(T^{n-1})}{T^n - T^{n-1}} + \frac{C^{n-1}(T^{n-1}) - C^n(T^{n-1})}{T^n - T^{n-1}} \nonumber \\
                    & \approx \frac{C^n(T^n) - C^{n-1}(T^{n-1})}{T^n - T^{n-1}}, \nonumber
\end{align}
where, in the last step, we have assumed that $C^n$ varies slowly as a function of $n$. This yields

\begin{equation}
    T^{n+1} = T^n - \gamma \frac{C^n(T^n) - C^{n-1}(T^{n-1})}{T^n - T^{n-1}}. \nonumber
\end{equation}
Taking exponentials on both sides of the equation, we get
\begin{equation}
    \Delta t^{n+1} = \Delta t^n \mathrm{exp}(-\gamma \Delta), \qquad \Delta = \frac{\mathrm{ln}\, c^n - \mathrm{ln}\, c^{n - 1}}{\mathrm{ln}\, \Delta t^n - \mathrm{ln}\, \Delta t^{n - 1}}. \nonumber
\end{equation}
Here, $\gamma$ is a free parameter, which can be a function of $c$ and $\Delta t$. Choosing $\gamma$ to be a constant, despite being the simplest choice, has two major disadvantages. First, we can not guarantee that $\Delta t^n \not = \Delta t^{n-1}$, which would lead to numerical problems in computing $\Delta$. Second, in some situations the controller can yield prohibitively large changes in step size. We, therefore, compute the new step size $\Delta t^{n+1}$ as follows
\[
  \Delta t^{n+1} = \Delta t^n \times
  \begin{cases}
  
    \lambda                                                 & \text{if $1 \leq s < \lambda$}, \\
    \delta                                                  & \text{if $ \delta \leq s < 1$}, \\
    s:=\exp({-\alpha \; \mathrm{tanh}(\beta \Delta)})       & \text{otherwise}.
  \end{cases}
\]
The parameter $\alpha$ acts as a constraint on the maximal allowed step-size change: the maximum change in step size is given by $\exp(\pm \alpha)$. The parameter $\beta$ determines how strongly the controller reacts to a change in the cost. The factors $\lambda$ and $\delta$ have been incorporated to ensure that the step size changes by at least $\lambda \, \Delta t$ or  $\delta \, \Delta t$ depending on whether $\Delta t$ needs to be increased or decreased for minimizing the cost. They are chosen in such a way that results in non-trivial changes in the step size if $\Delta$ is close to 0. 

These parameters, ($\alpha$, $\beta$, $\lambda$, and $\delta$) have been numerically optimized for the linear diffusion-advection equation using an implicit Runge--Kutta scheme for a range of values of $N$, $\eta$, and $\text{tol}$. This step-size controller has been designed with two variants: (i) \textit{Non-penalized} variant: the aforementioned parameters have been chosen to incur the minimum possible cost, (ii) \textit{Penalized} variant: if the traditional controller performs better than the non-penalized variant, a penalty is imposed on the proposed controller. This penalty has been imposed to trade off the enhanced performance of the proposed controller (where it performs better) with an acceptable amount of diminished performance (where its performance is inferior to that of the traditional controller). The numerical optimization yields the following set of parameters
\begin{align*}
    \textbf{Non-penalized} \quad \alpha = 0.65241444 \quad \beta = 0.26862269 \quad \lambda = 1.37412002 \quad \delta = 0.64446017 \\
    \textbf{Penalized} \quad \alpha = 1.19735982  \quad \beta = 0.44611854 \quad \lambda = 1.38440318 \quad \delta = 0.73715227
\end{align*}
The improved performance of this controller for both implicit \cite{Einkemmer18} and exponential integrators (as we will see in this paper) further shows the generality of this approach. This is emphasized as the parameters that have been obtained for a linear PDE generalize well to nonlinear problems for a variety of numerical methods.

The different strategy in these two approaches can be understood from the coefficients of the non-penalized and penalized controller. The penalized variant tends to change the step size more rapidly. In cases where the step size is limited by accuracy considerations, less exploration (optimization) is performed, and it is expected to behave almost identical to the classic step size controller. The non-penalized step-size controller, however, explores more possibilities and is generally (as we will see) more efficient in the majority of cases (where the step size is chosen due to considerations of computational cost). Details about the optimization procedure can be found in \cite{Einkemmer18} (cf. Fig.~1).

It is worth noting that the proposed step-size controller is designed solely to minimize the computational cost and does not take into account the error incurred during the time integration. As such, we need to take the minimum of the step sizes given by the traditional and the proposed step-size controller
\begin{equation}
    \Delta t^{n + 1} = \mathrm{min} \left(\Delta t^{n + 1}_\mathrm{proposed}, \; \Delta t^{n + 1}_\mathrm{traditional}\right) \nonumber
\end{equation}
This ensures that, in addition to satisfying the accuracy requirements set by the user, the controller minimizes the computational cost (by choosing a smaller time step size) whenever possible.

In the context of exponential integrators, there is an additional detail that needs to be taken care of: the computation of the Leja interpolation to the prescribed tolerance. In some situations, especially if enormously large step sizes are chosen, the polynomial interpolation might fail to converge within a reasonable number of iterations. If this is the case for any of the internal stages of the exponential Rosenbrock scheme, we reject the step and use the traditional controller to determine a smaller step size.

%%%%%%%%%%%%%%%%%%%%%%%%%%%%%%%%%%%%%%%%%%%%%%%%%%%%%%%%%%%%%%%%%%%%%%%%%%%%%%%%%%%%%%%%%%%%%%%%%%%%%%%%%%%

\section{Numerical results}
\label{sec:nl_eqs}

In this section, we investigate the performance of the proposed step-size controller and compare it with the traditional controller using the embedded EXPRB43 scheme for a number of nonlinear problems. We present a detailed explanation of how the step-size controller can improve the performance of exponential Rosenbrock integrators. Most of the examples in this section have been adopted from \cite{Einkemmer18}. In all of these examples, we consider periodic boundary conditions on $[0, 1]$ and $[0, 1] \times [0, 1]$ for the 1D and 2D cases respectively.

%%% ----------------------------------------------------------- %%%
%%% ----------------------------------------------------------- %%%

\subsection{Viscous Burgers' Equation}

\begin{figure*}[htp]
    \centering
	\includegraphics[width = \textwidth]{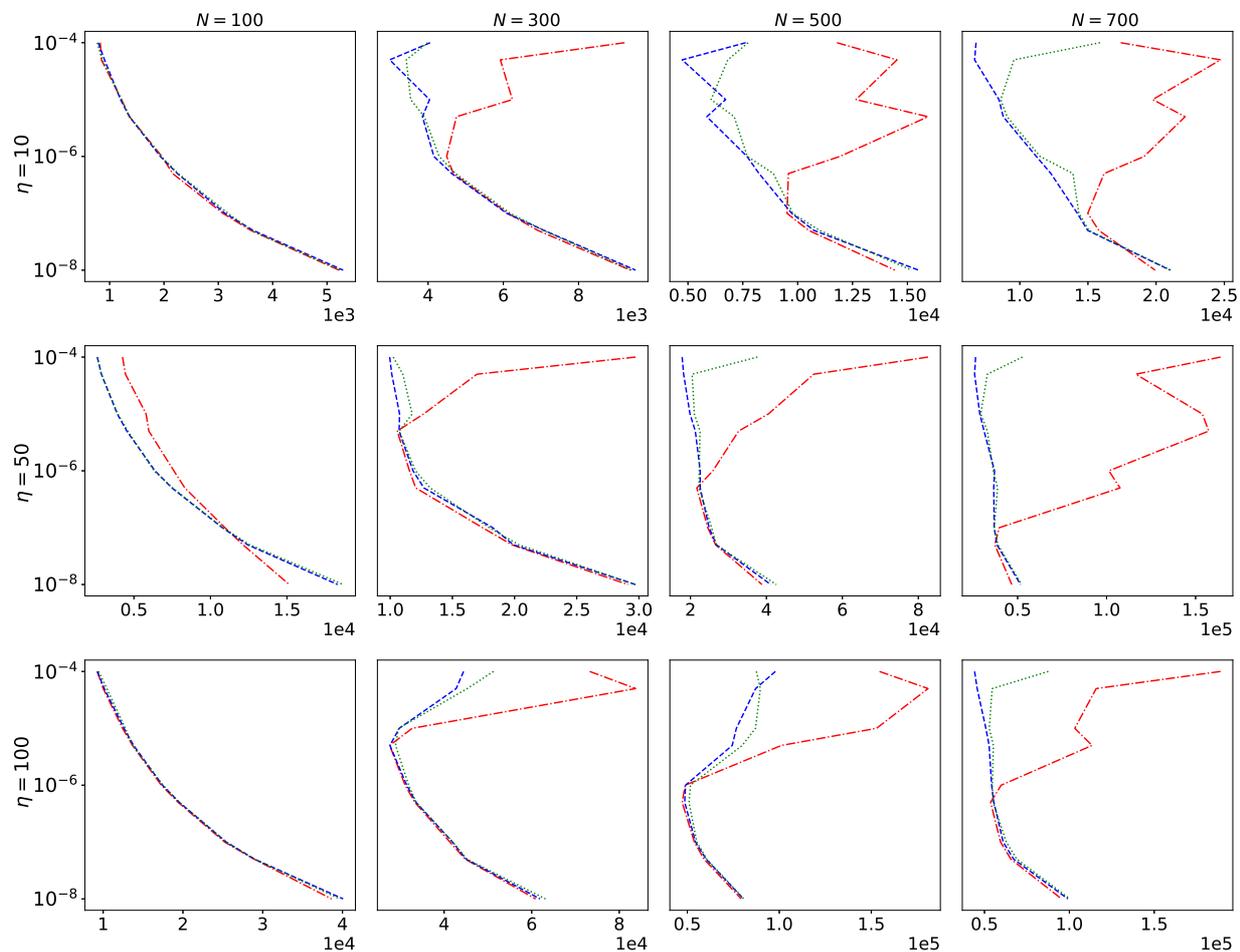}
    \caption{Figure shows the number of matrix-vector products (i.e. computational cost) vs. tolerance for the 1D viscous Burgers' equation for different values of $\eta$ and $N$. The red dashed-dotted lines correspond to the traditional controller, green dotted lines the penalized variant, and the blue dashed lines represent the non-penalized variant.}
    \label{fig:exprb43_viscous_1d}
\end{figure*}

\begin{figure*}
    \centering
	\includegraphics[width = 0.8\textwidth]{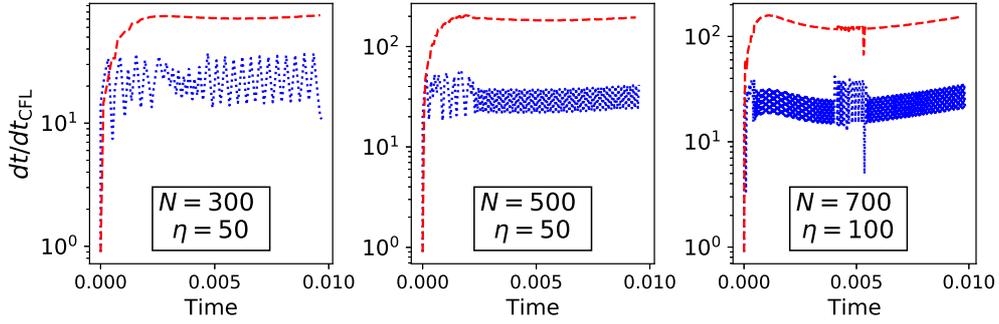}
    \caption{We compare the step sizes used during the simulation with the non-penalized variant of the proposed controller (blue dotted lines) with the step sizes yielded by the traditional controller for tolerance, $\text{tol}=10^{-4}$. The step sizes are significantly reduced which results in improved performance. Here, $dt_\mathrm{CFL} = \mathrm{min}(1/2N^2, 1/\eta\,N)$.}
    \label{fig:dt_compare_used_trad}
\end{figure*}

\begin{figure*}
    \centering
	\includegraphics[width = 0.7\textwidth]{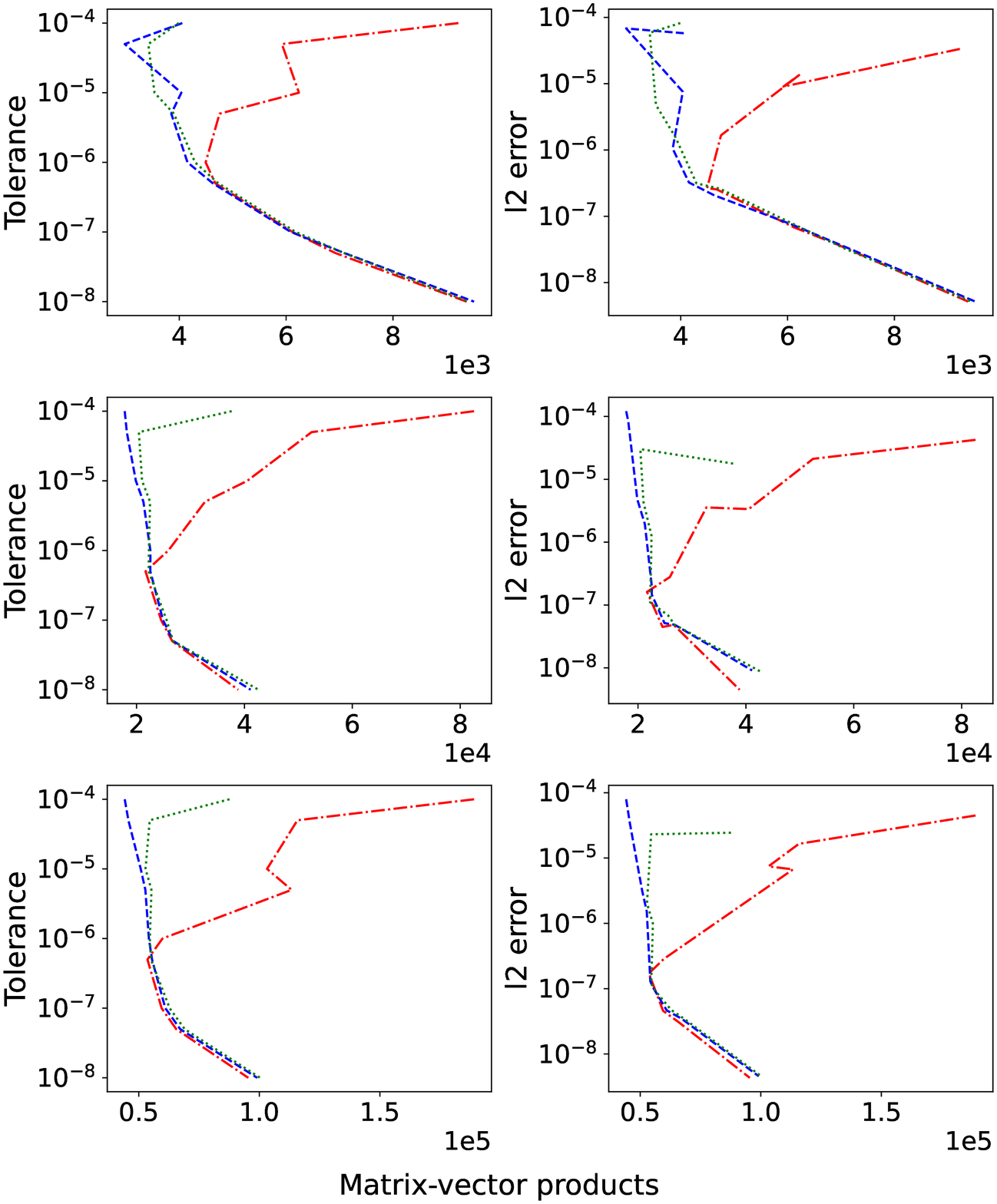}
    \caption{A comparison of the the l2 norm of the error incurred (right column) as a function of computational cost (i.e. number of matrix-vector products) for given values of user-defined tolerances (left column) for 3 different cases: (i) $N = 300, \eta = 10$ (top panel), (ii) $N = 500, \eta = 50$ (middle panel), and (iii) $N = 700, \eta = 100$ (bottom panel) for the 1D viscous Burgers' equation. The red dashed-dotted lines correspond to the traditional controller, green dotted lines the penalized variant, and the blue dashed lines represent the non-penalized variant. We can clearly see that the accuracy constraints are well satisfied by all the step-size controllers.}
    \label{fig:exprb43_viscous_1d_error}
\end{figure*}

The one-dimensional viscous Burgers' equation (conservative form) reads
\begin{equation}
    \frac{\partial u}{\partial t} = \frac{1}{2} \eta \frac{\partial u^2}{\partial x} + \frac{\partial^2 u}{\partial x^2}, \nonumber
\end{equation}
The P\'eclet number ($\eta$) is a measure of the relative strength of advection to diffusion. Higher values of $\eta$ indicate advection-dominated scenarios whereas lower values of $\eta$ imply diffusion-dominated cases. The initial condition is 
\begin{equation}
	u(x, t = 0) = 1 + \mathrm{exp}\left(1 - \frac{1}{1-(2x - 1)^2}\right) + \frac{1}{2} \mathrm{exp}\left(-\frac{(x - x_0)^2}{2\sigma^2}\right)  \nonumber
\end{equation}
with $x_0 = 0.9$ and $\sigma = 0.02$. 

We wish to test our step-size controller in diffusion as well as advection dominated cases. If $\eta$ is small (diffusion dominated), then the Gaussian part is dynamically smeared out, and $u(x, t)$ is slowly advected over a significant amount of time.  If $\eta$ is large (advection dominated), the solution undergoes rapid advection and the Gaussian is smeared out after a long time. We choose the final time of the integration long enough ($t = 10^{-2}$) such that for any value of $\eta $, a fixed amount
of diffusion is inherently introduced in the simulations. For the space discretization, we consider a third-order upwind scheme for advection and the second-order centered difference scheme for diffusion (see Appendix \ref{app:diff_adv}).

The work-precision diagram is shown in Fig. \ref{fig:exprb43_viscous_1d}. The computational cost incurred, in our case, measured by the number of matrix-vector products \footnote{One can also choose the number of Leja points used at every time step as a proxy of the computational cost}, is plotted as a function of the user-specified tolerance. The blue curves correspond to the non-penalized variant of the proposed controller, the green ones refer to the penalized variant, and the red lines represent the traditional controller. Let us clearly state the non-penalized variant is the recommended step-size controller. The idea of the penalized variant is to have a controller that behaves similar to the traditional one in cases where the classic approach is advantageous. However, as our numerical simulations show, this is almost never the case. Overall, the non-penalized controller has significantly enhanced performance.

The performance of the proposed step-size controller is similar to the traditional controller for $N = 100$, where $N$ is the number of grid points, for all values of $\eta$. As the number of grid points is increased, one can see that the proposed step-size controller performs significantly better in the lenient to intermediate tolerance regime. This is true for all values of $\eta$. Maximum speedups up to a factor of 2.5 are observed. To illustrate how this performance improvement is achieved, we compare, in Fig. \ref{fig:dt_compare_used_trad}, the step size used at each time step during the simulations with the proposed controller (blue curves) with the largest possible step size (constrained only by the accuracy requirements) estimated by the traditional controller (red curves). We see that the step sizes estimated by the proposed controller are smaller than what would be possible based purely on accuracy constraints. This justifies the fundamental principle of the proposed step-size controller; i.e., multiple small step sizes incur less computational effort than a single large step size. It can also be seen in Fig. \ref{fig:dt_compare_used_trad} that the proposed step-size controller continuously varies the step size to find the step size that minimizes the computational cost. For stringent tolerances, the traditional controller already yields small step sizes. Further reduction in step sizes would only result in an increased number of time steps leading to an increase in the computational cost. All in all, the step-size controller (both non-penalized and penalized variants) has superior performance for a reasonably wide range of tolerances for advection as well as diffusion dominated cases. Moreover, as the step sizes are further reduced by the proposed controller, the solution is, in fact, more accurate. This is in addition to the reduced computational cost. It also largely avoids the inverse C-shaped curve observed for the traditional controller. We note that the proposed step-size controller is notably efficient in the lenient to intermediate tolerance regime which is relevant for most practical applications. To further validate our simulation results, we show, in Fig.~\ref{fig:exprb43_viscous_1d_error}, the error incurred as a function of the computational cost for three different cases. It can be seen that the error incurred remains at par, if not below, the user-defined tolerances. One can also appreciate the difference in the computational times between the proposed controller and the traditional controller in the lenient to intermediate tolerance range.

\begin{figure*}
    \includegraphics[width = \textwidth]{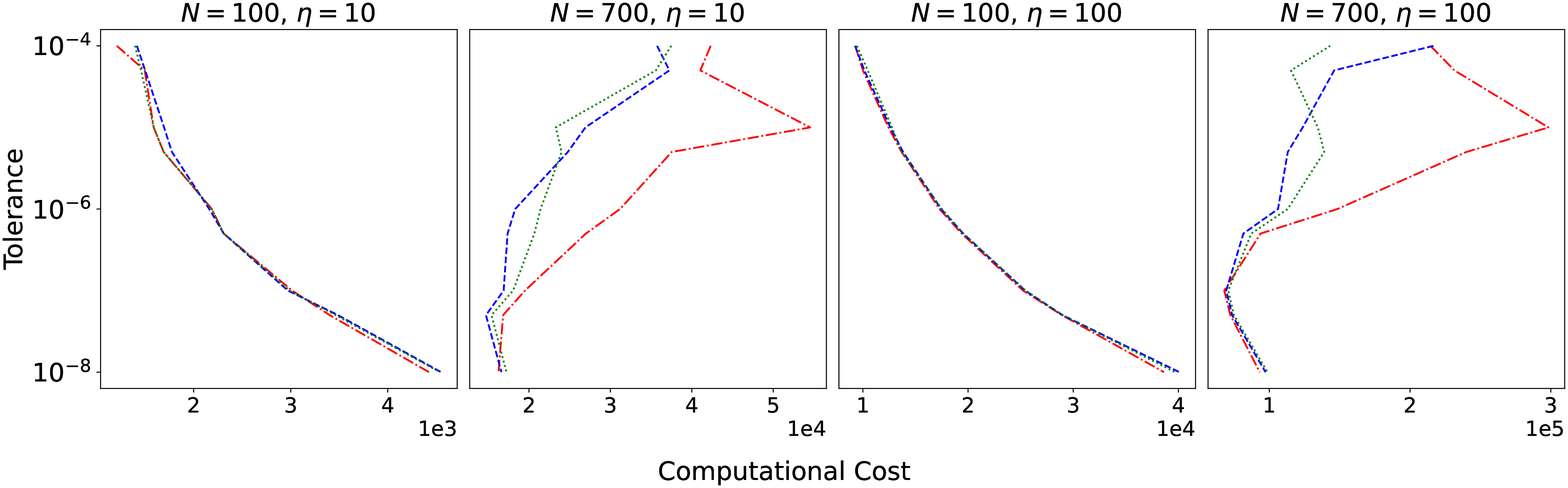}
    \includegraphics[width = \textwidth]{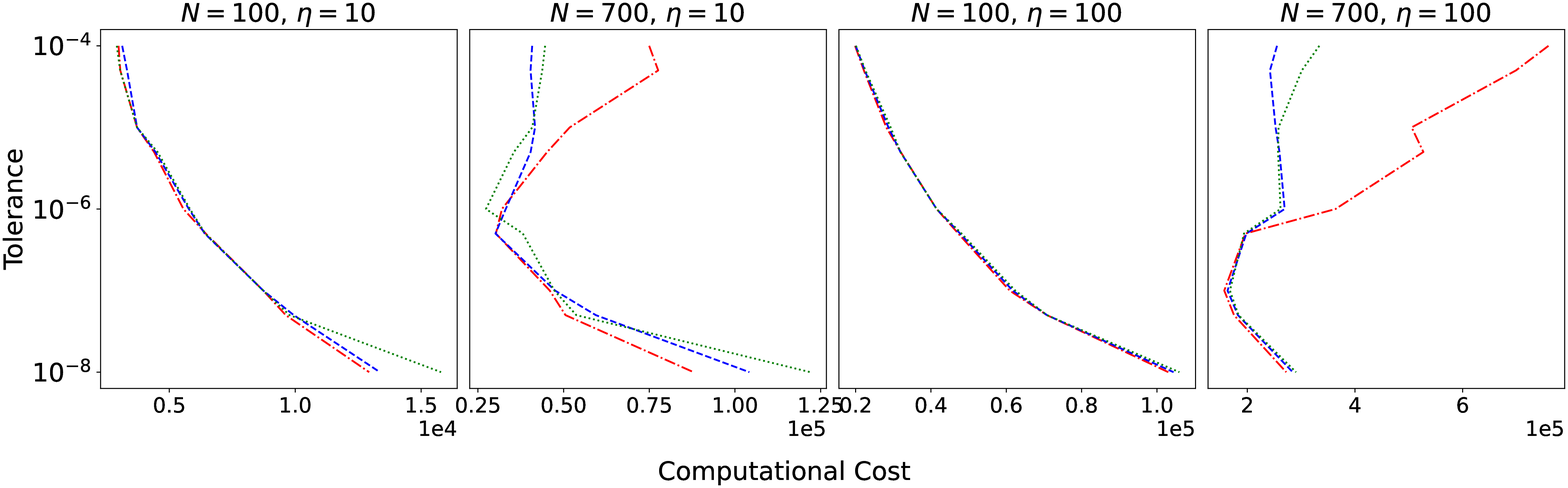}
    \includegraphics[width = \textwidth]{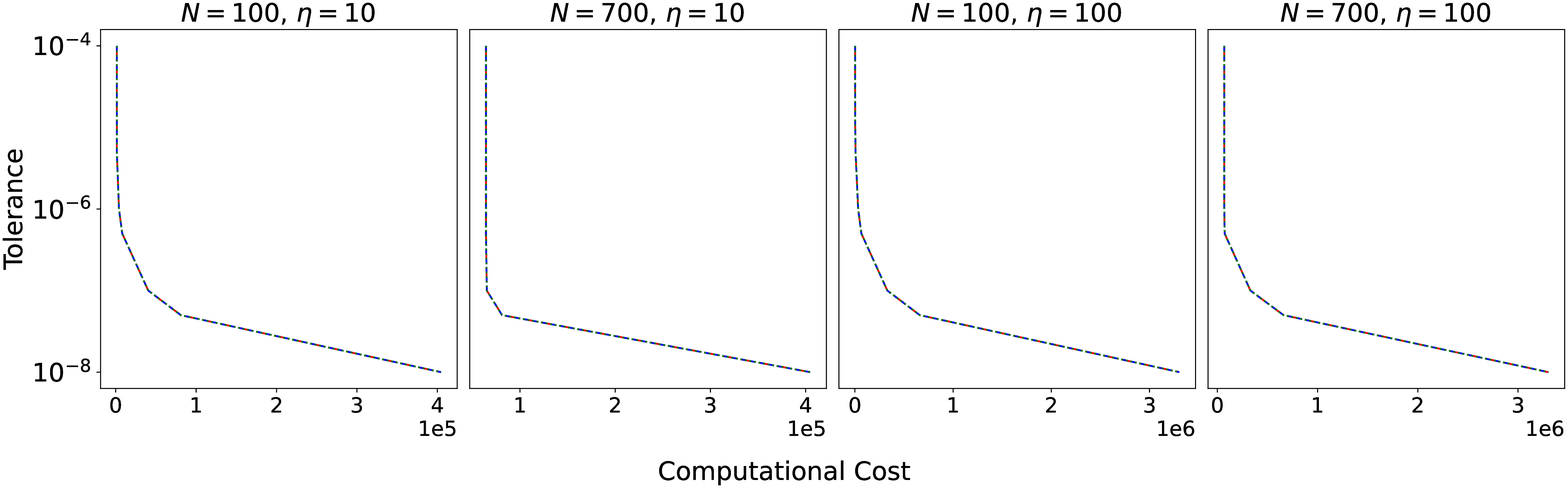}
    \caption{Comparison of the number of matrix-vector products (i.e. computational cost) vs. tolerance for the 1D viscous Burgers' equation for the embedded EXPRB43 scheme (top panel), Richardson extrapolation with the third-order solution of this embedded scheme (middle panel), and the explicit RKF45 integrator (bottom panel). Here we consider two different values of $\eta$ and $N$ each. The red dashed-dotted lines correspond to the traditional controller, green dotted lines the penalized variant, and the blue dashed lines represent the non-penalized variant.}
    \label{fig:compare_RE_embed_RKF}
\end{figure*}

Using an embedded method is not the only way to obtain an error estimate for automatic step size control. Richardson extrapolation, whilst usually being more expensive, has the advantage that it can be applied to any time integration scheme. The proposed step-size controller is independent of how the error estimate is obtained. To illustrate this, we apply it to the third-order solution of EXPRB43 using Richardson extrapolation as an error estimator. The results are presented in Fig. \ref{fig:compare_RE_embed_RKF}, along with a comparison with the embedded EXPRB43 and the explicit embedded Runge--Kutta--Fehlberg 45 (RKF45) schemes. It can clearly be seen that the proposed step-size controller works reasonably well whilst using Richardson extrapolation as an error estimator. However, the embedded Rosenbrock scheme preemptively outperforms the Richardson extrapolation method, as is expected. It is worth noting that the `shape' of the curves is fairly similar for both these methods. This tells us that the embedded scheme and the Richardson extrapolation have similar changes in behaviour with the increase or decrease in tolerance. The explicit embedded scheme RKF45 (fourth-order error estimate) is over an order of magnitude more expensive than the corresponding exponential integrator counterparts. Any reduction in step size, over the ones given by the traditional controller, would only result in an increased number of time steps leading to an increase in the computational cost. As such, the performance of the traditional controller and the proposed controller is the same for this integrator. 

\begin{figure}[t]
    \centering
	\includegraphics[width = 0.8\textwidth]{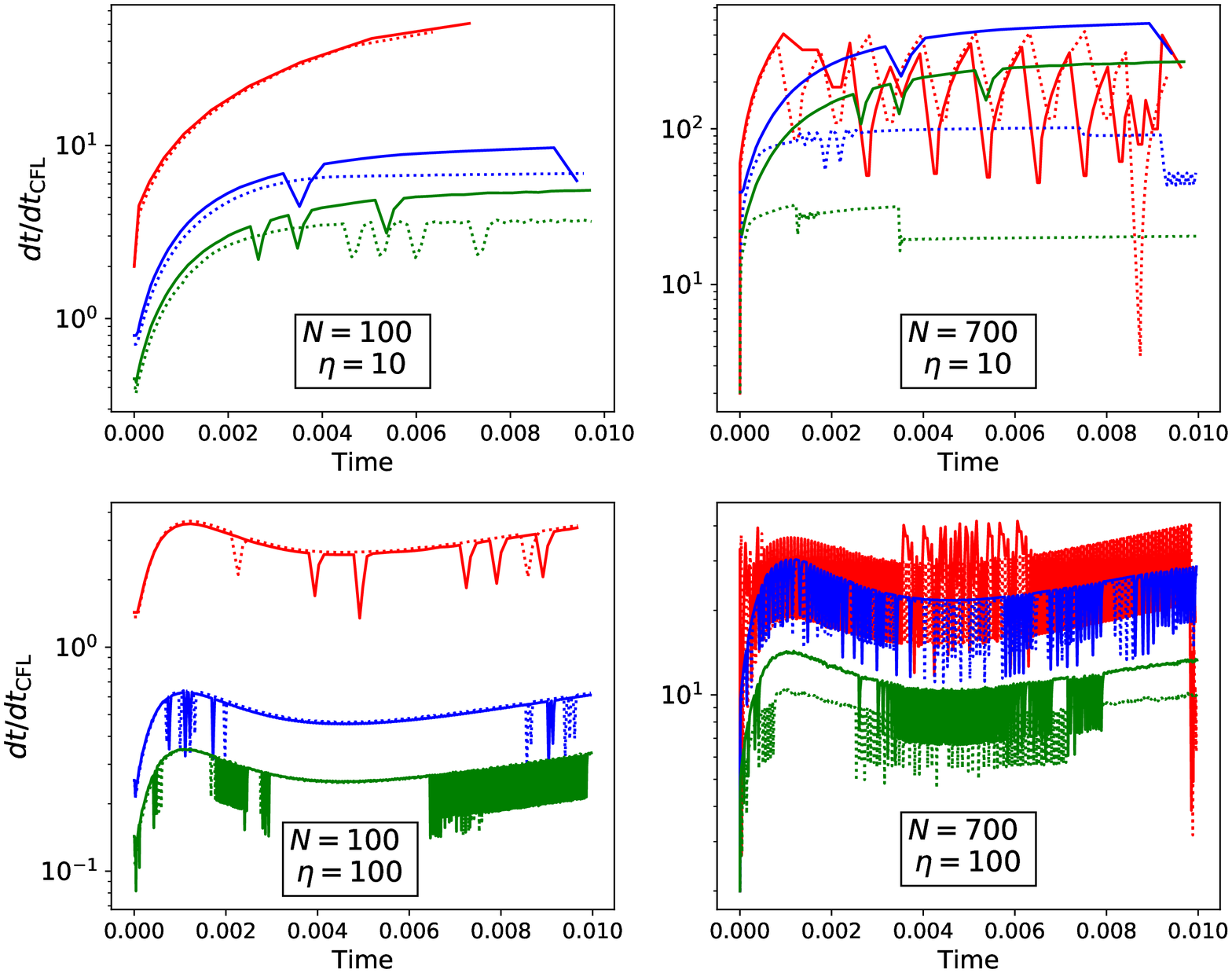}
    \caption{Figure illustrates the step sizes (normalized to the CFL time) for four different combinations of $N$ and $\eta$. The solid lines represent the embedded EXPRB43 scheme and the dotted lines represent Richardson extrapolation with the third-order solution. The colours indicate different tolerances: $10^{-4}$ (red), $10^{-7}$ (blue), and $10^{-8}$ (green). The step sizes are shown for the non-penalized variant.}
    \label{fig:dt_compare}
\end{figure}

\begin{figure*}
    \centering
	\includegraphics[width = 0.8\textwidth]{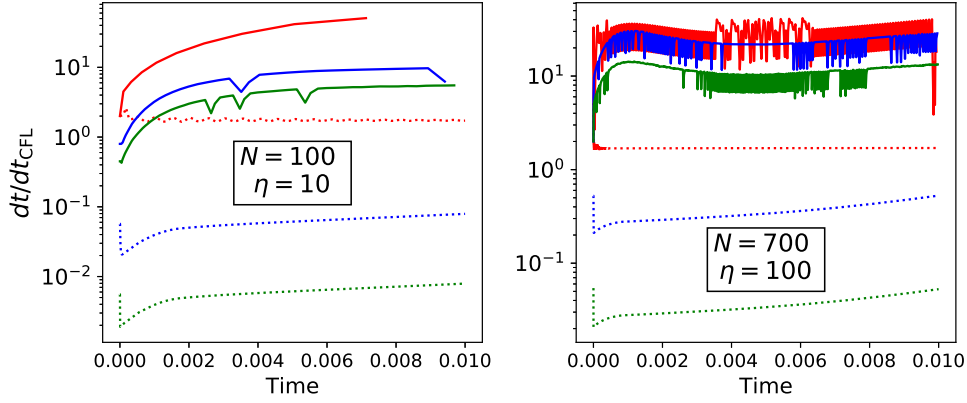}
    \caption{Figure illustrates the step sizes (normalized to the CFL time) for two different combinations of $N$ and $\eta$. The solid lines represent the embedded EXPRB43 scheme and the dotted lines represent the explicit RKF45. The colours indicate different tolerances: $10^{-4}$ (red), $10^{-7}$ (blue), and $10^{-8}$ (green). The step sizes are shown for the non-penalized variant.}
    \label{fig:dt_compare_rkf}
\end{figure*}

Fig. \ref{fig:dt_compare} compares the step sizes (non-penalized variant) for the embedded EXPRB43 scheme and Richardson extrapolation with the third-order solution. One can see that the step sizes for the embedded scheme are, in general, larger than the step sizes for the Richardson extrapolation. The increased number of time steps is likely to be an additional contribution to the expenses of the Richardson extrapolation method. For $N = 700; \eta = 10$, the step sizes for the embedded scheme are similar to that of Richardson extrapolation in the lenient tolerance regime ($\mathrm{tol} = 10^{-4}$: red lines). As such, the computational costs are somewhat similar. As the tolerance is reduced, the step sizes allowed by Richardson extrapolation decrease substantially, thereby incurring more cost ($\mathrm{tol} = 10^{-7} \; \text{and} \; 10^{-8}$: blue and green lines respectively). Similar arguments can be used to explain the other cases as well. The step sizes permitted by the explicit RKF45, depicted in Fig. \ref{fig:dt_compare_rkf}, are significantly smaller (roughly 1 - 2 orders of magnitude) than EXPRB43. Consequently, this incurs a hefty computational cost and has the worst performance out of the three schemes presented here.

\begin{table}
\centering
\subfloat[]{
\begin{tabular}{l|c|c|c}
\textbf{Parameters} & \textbf{SDIRK23} & \textbf{EXPRB43} & \textbf{RKF45}\\ 
\hline
    $N = 100, \eta = 10$    & $10^4 - 3 \cdot 10^4$             & $10^3 - 4 \cdot 10^3$              & $1.5 \cdot 10^3 - 4 \cdot 10^5$ \\ 
	$N = 100, \eta = 100$ 	& $5 \cdot 10^4 - 2 \cdot 10^5$     & $2 \cdot 10^4 - 3 \cdot 10^4$      & $6 \cdot 10^4 - 4 \cdot 10^5$ \\
	\hline
	$N = 700, \eta = 10$    & $5 \cdot 10^4 - 1.5 \cdot 10^5$   & $10^4 - 4 \cdot 10^4$              & $2 \cdot 10^3 - 3.5 \cdot 10^6$ \\ 
    $N = 700, \eta = 100$   & $5 \cdot 10^5 - 1.5 \cdot 10^6$   & $10^5 - 2 \cdot 10^5$              & $7 \cdot 10^4 - 3.5 \cdot 10^6$ \\
\end{tabular}}
\caption{A quantitative comparison of the computational cost incurred by EXPRB43, RKF45, and SDIRK23 (used in \cite{Einkemmer18}) for the viscous Burgers' equation within the range of tolerance ($10^{-4} - 10^{-8}$) considered in this study.}
\label{tab1}
\end{table}

Now, we compare the performance of the exponential Rosenbrock approach with implicit and explicit integrators. Table \ref{tab1} compares the performance of the non-penalized variant for the embedded EXPRB43 scheme with the two-stage third-order singly diagonally implicit Runge-Kutta (SDIRK23) scheme used in \cite{Einkemmer18} (with the same controller) and the RKF45 scheme. It is evident that the embedded Rosenbrock method has superior performance compared to the other two; up to an order of magnitude over the implicit integrator and up to two orders of magnitude over the explicit integrator can be observed for some configurations. Similar results have been obtained for the inviscid Burgers' equation and the porous medium equation that are discussed in the following sections.

%%% ----------------------------------------------------------- %%%

Next, we test our step-size controller on the two-dimensional viscous Burgers' equation:
\begin{equation}
    \frac{\partial u}{\partial t} = \frac{1}{2} \left(\eta_\mathrm{x} \frac{\partial u^2}{\partial x} + \eta_\mathrm{y} \frac{\partial u^2}{\partial y}\right) + \frac{\partial^2 u}{\partial x^2} + \frac{\partial^2 u}{\partial y^2}, \nonumber
\end{equation}
where $\eta_\mathrm{x}$ and $\eta_\mathrm{y}$ are the components of the P\'eclet number along the $X$ and $Y$ directions, respectively. The initial condition is chosen to be

\begin{equation}
	u(x, y, t = 0) = 1 + \mathrm{exp}\left(1 - \frac{1}{1-(2x - 1)^2} - \frac{1}{1-(2y - 1)^2}\right) + \frac{1}{2} \mathrm{exp}\left(\frac{-(x - x_0)^2 - (y - y_0)^2}{2\sigma^2}\right)  \nonumber
\end{equation}
with $x_0 = 0.9$, $y_0 = 0.9$, and $\sigma = 0.02$.

\begin{figure*}
    \centering
	\includegraphics[width = 0.975\textwidth]{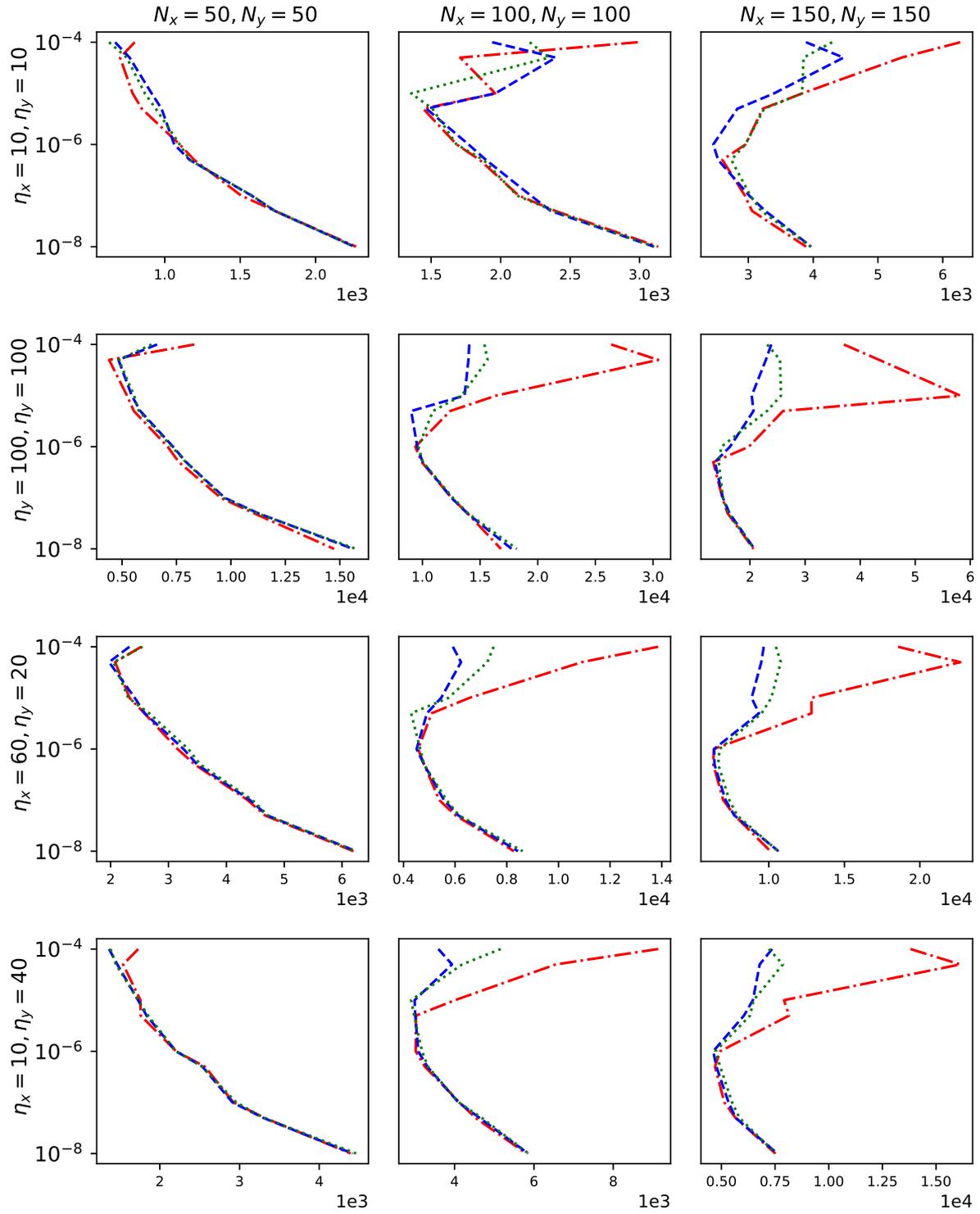}
    \caption{Figure shows the number of matrix-vector products (i.e computational cost) vs. tolerance for the 2D viscous Burgers' equation. The values of $N_x$, $N_y$, $\eta_x$, and $\eta_y$ are varied. The red dashed-dotted lines correspond to the traditional controller, green dotted lines the penalized variant, and the blue dashed lines represents the non-penalized variant of the proposed step-size controller.}
    \label{fig:exprb43_viscous_2d}
\end{figure*}

The work-precision diagram is shown in Fig. \ref{fig:exprb43_viscous_2d}. The proposed step-size controller has a similar performance compared to the one-dimensional case. With the increase in the number of grid points ($N_x$ and $N_y$ correspond to the number of grid points along $X$ and $Y$ directions respectively), the proposed controller shows a large improvement in performance (up to a factor of 3) for lenient tolerances for the different values of $\eta_x$ and $\eta_y$ considered here. 

%%% ----------------------------------------------------------- %%%
%%% ----------------------------------------------------------- %%%

\subsection{Inviscid Burgers' Equation}

\begin{figure}
    \centering
	\includegraphics[width = 0.9\textwidth]{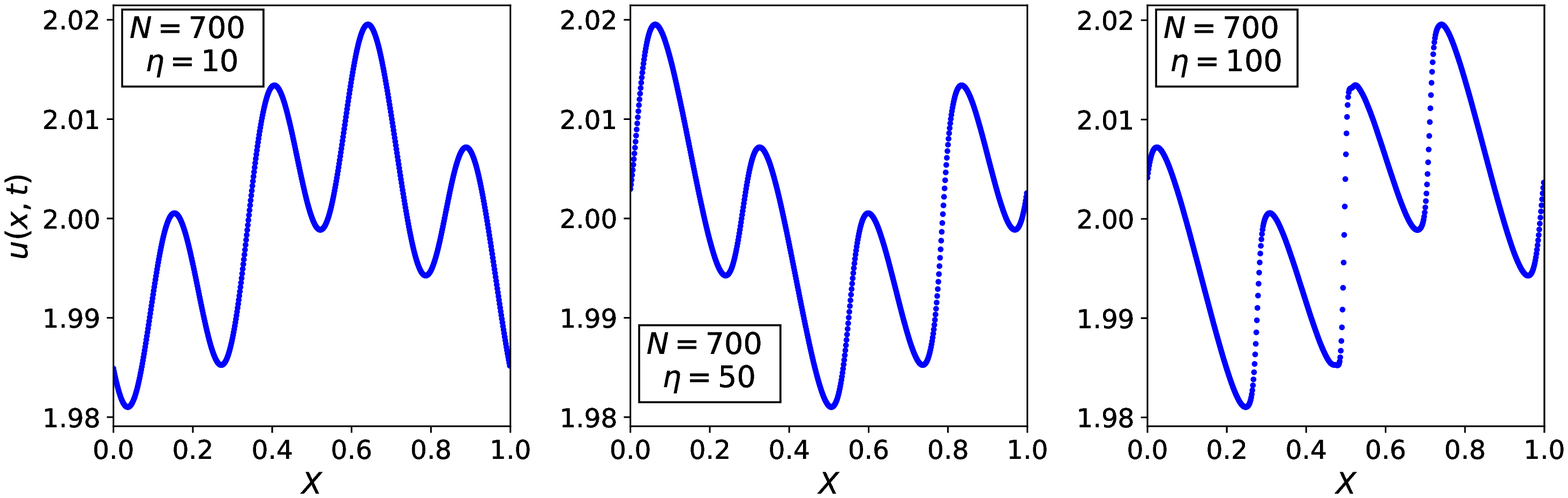}
    \caption{Figure shows the solution at the final time $t = 3.25 \eta \times 10^{-2}$ for the inviscid Burgers' equation. It can be seen that the solution yields progressively steeper gradients for large values of $\eta$ , i.e. the P\'eclet number. Further increase in the simulation time would result in a shock.}
    \label{fig:final_u_invis_1d}
\end{figure}

We consider the conservative form of the inviscid Burgers' equation
\begin{equation}
    \frac{\partial u}{\partial t} = \frac{1}{2} \frac{\partial u^2}{\partial x}. \nonumber
\end{equation}
The initial condition is given by
\begin{equation}
    u(x, t = 0) = 2 + \epsilon_1 \sin(\omega_1 x) + \epsilon_2 \sin(\omega_2 x + \varphi) \nonumber
\end{equation}
with $\epsilon_1 = \epsilon_2 = 10^{-2}$, $\omega_1 = 2\pi$, $\omega_2 = 8\pi$, and $\phi = 0.3$. The simulations are carried out until $t = 3.25\eta \times 10^{-2}$, where $\eta$ is the P\'eclet number. The distribution at the final time, for different values of $\eta$, is depicted in Fig. \ref{fig:final_u_invis_1d}. A change in $\eta$ corresponds to a change in the final time of the simulation. As time progresses, the gradients start becoming progressively sharper. The solution gradually approaches a shock wave. It is worth noting that this effect is more prominent in cases with larger values of $N$. This is due to the fact that the numerical diffusion decreases as the number of grid points are increased leading to increasingly steeper gradients.

\begin{figure}
    \centering
	\includegraphics[width = 0.9\textwidth]{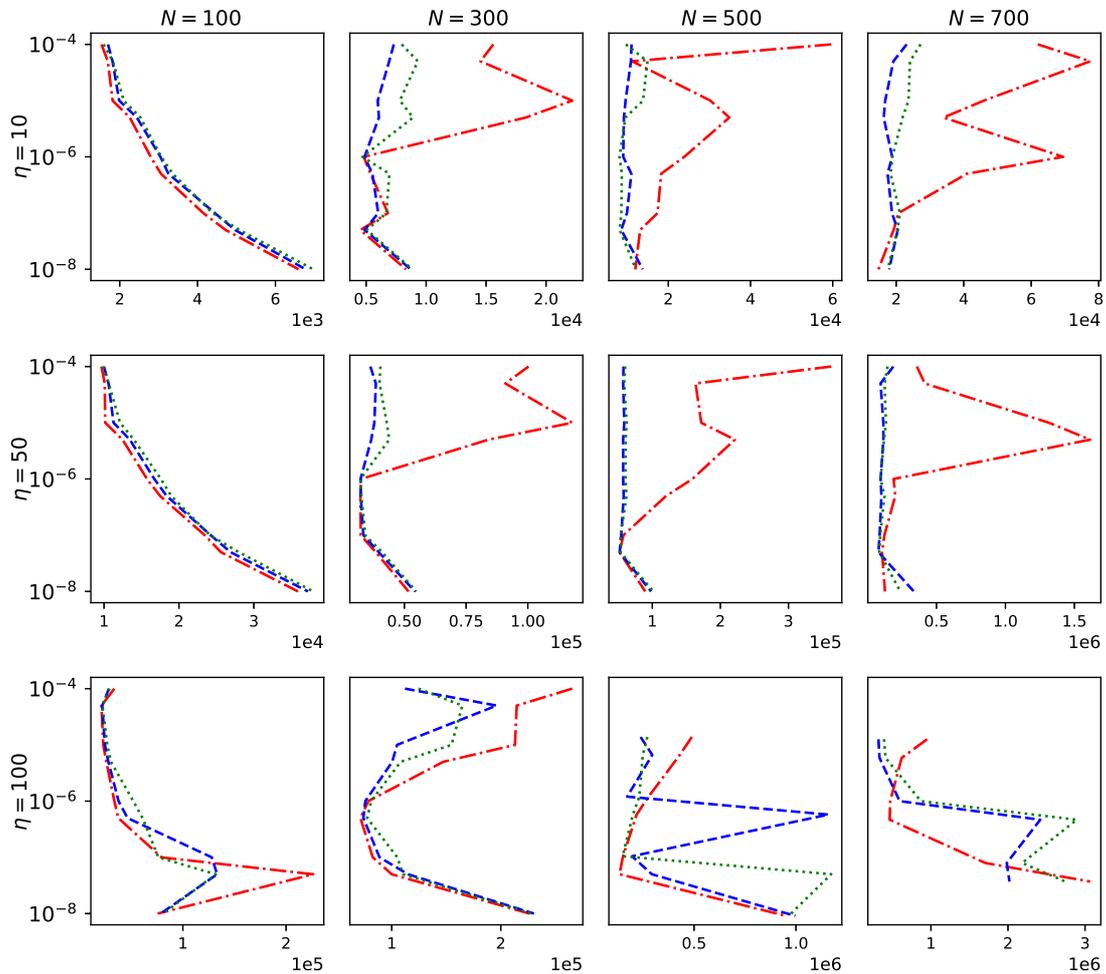}
    \caption{Figure shows the number of matrix-vector products (i.e computational cost) vs. tolerance for the 1D inviscid Burgers' equation for different values of $\eta$ and $N$. The red dashed-dotted lines correspond to the traditional controller, green dotted lines the penalized variant, and the blue dashed lines represents the non-penalized variant of the proposed step-size controller.}
    \label{fig:exprb43_inviscid_1d}
\end{figure}

The work-precision diagram is illustrated in Fig. \ref{fig:exprb43_inviscid_1d}. It can be seen that the proposed step-size controller `flattens out' the zig-zag shape of the curves yielded by the traditional controller to a large extent. This yields a significant improvement over the traditional controller,  especially in the lenient to medium tolerance range. Performance improvements of up to a factor of $4$ are observed. Both step-size controllers have some difficulty dealing with large values of $\eta$ and large $N$, i.e.~with very sharp gradients in the solution.

\begin{figure*}
    \centering
	\includegraphics[width = \textwidth]{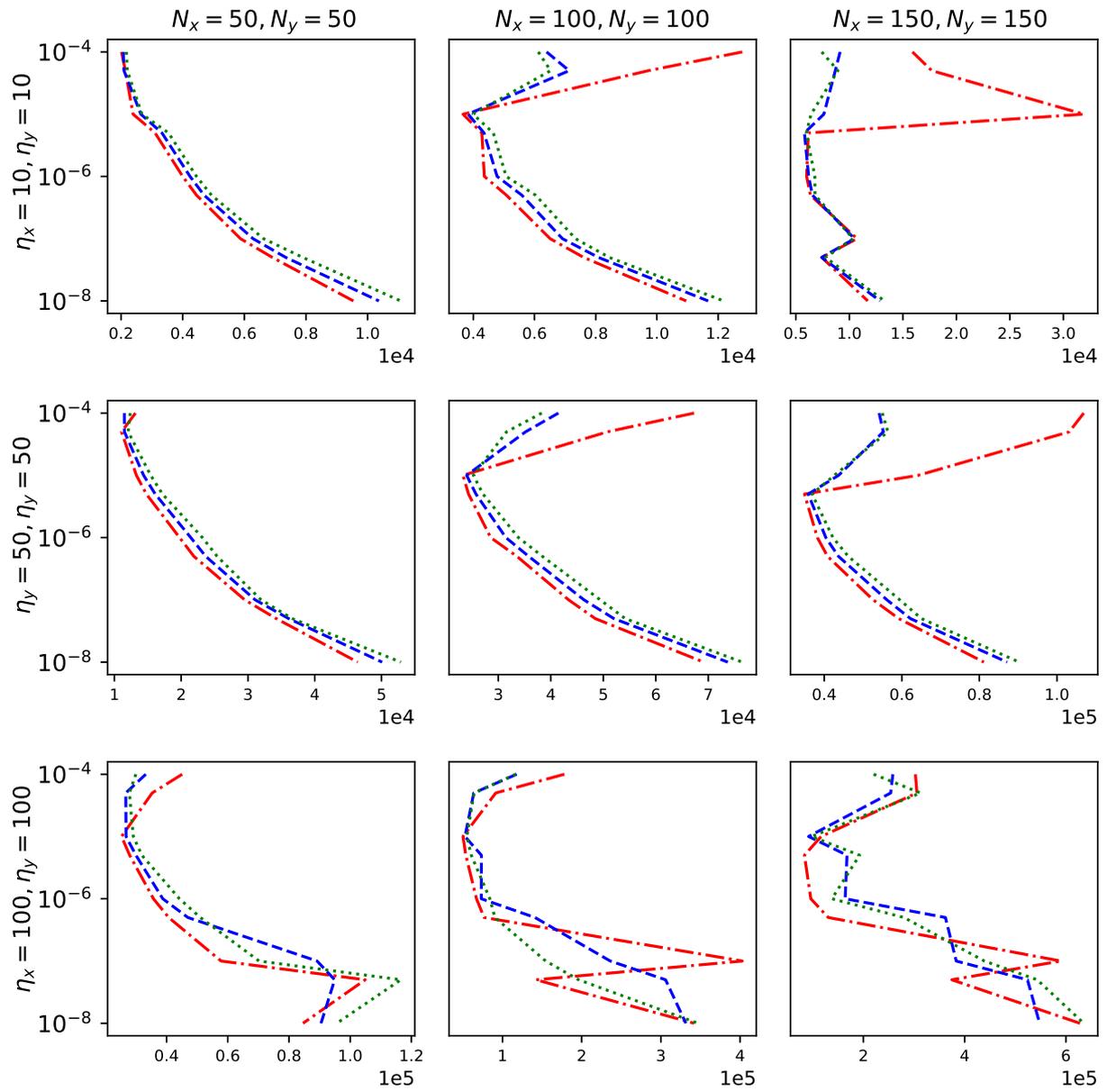}
    \caption{Figure shows the number of matrix-vector products (i.e computational cost) vs. tolerance for the 2D inviscid Burgers' equation. The values of $N_x$, $N_y$, $\eta_x$, and $\eta_y$ are varied. The red dashed-dotted lines correspond to the traditional controller, green dotted lines the penalized variant, and the blue dashed lines represents the non-penalized variant.}
    \label{fig:exprb43_inviscid_2d}
\end{figure*}

%%% ----------------------------------------------------------- %%%

We extend our 1D model into 2 dimensions. The two-dimensional inviscid Burgers' equation is given by
\begin{equation}
    \frac{\partial u}{\partial t} = \frac{1}{2} \left(\frac{\partial u^2}{\partial x} + \frac{\partial u^2}{\partial y}\right). \nonumber
\end{equation}
where we consider the initial condition
\begin{equation}
    u(x, y, t = 0) = 2 + \epsilon_1 \sin(\omega_1 x) + \epsilon_2 \sin(\omega_2 x + \phi) + \epsilon_1 \sin(\omega_1 y) + \epsilon_2 \sin(\omega_2 y + \phi) \nonumber
\end{equation}
with $\epsilon_1 = \epsilon_2 = 10^{-2}$, $\omega_1 = 2\pi$, $\omega_2 = 8\pi$, and $\phi = 0.3$. 

The corresponding work-precision diagram is shown in Fig \ref{fig:exprb43_inviscid_2d}. Once again, we see features similar to the one-dimensional case. The curves are flattened-out for lenient tolerances signifying a significant improvement over the traditional controller. For stringent tolerances, both controllers work well. It can also be seen that an increase in the number of grid points correlates with an enhanced performance of the proposed controller for a wide range of tolerance. This is also in agreement with what we have seen for the 1D case.

%%% ----------------------------------------------------------- %%%
%%% ----------------------------------------------------------- %%%

\subsection{Porous Medium Equation}

\begin{figure*}[t]
    \centering
	\includegraphics[width = 0.9\textwidth]{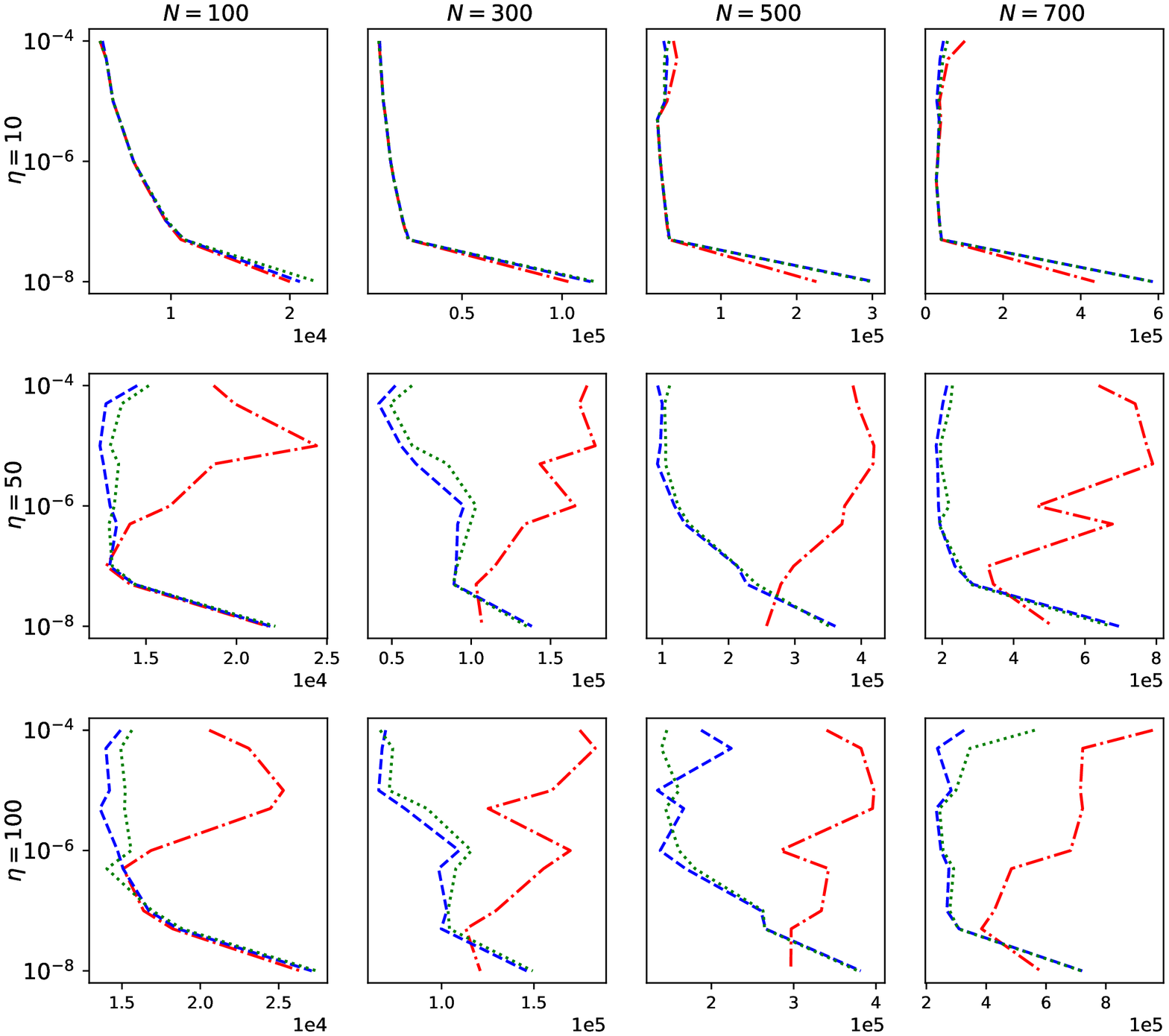}
    \caption{Figure shows the number of matrix-vector products (i.e. computational cost) vs. tolerance for the 1D porous medium equation for different values of $\eta$ and $N$. The red dashed-dotted lines correspond to the traditional controller, green dotted lines the penalized variant, and the blue dashed lines represent the non-penalized variant of the proposed step-size controller.}
    \label{fig:exprb43_porous_1d}
\end{figure*}

The next example considered is the porous medium equation with linear advection. The one-dimensional equation reads
\begin{equation}
    \frac{\partial u}{\partial t} = \eta \frac{\partial u}{\partial x} + \frac{\partial^2 u^m}{\partial x^2}, \nonumber
\end{equation}
where we have chosen $m = 2$ and $\eta$ is the P\'eclet number. Here, the initial condition is given by
\begin{equation}
    u(x, t = 0) = 1 + \Theta(x_1 - x) + \Theta(x - x_2). \nonumber
\end{equation}
Here, $x_1 = 0.25$, $x_2 = 0.6$, and $\Theta$ is the Heaviside function. This corresponds to a rectangle, i.e.~a discontinuous initial value. The nonlinear diffusivity dynamically smears out this discontinuity as the system evolves in time and results in a smooth solution. The simulations are carried out up to a final time of $t = 10^{-2}$.

The corresponding results are shown in Fig. \ref{fig:exprb43_porous_1d}. The performance of the proposed controller is similar to that of the traditional controller for $\eta = 10$ in the lenient to medium tolerance range. As $\eta$ increased, one can appreciate the significant reduction in computational cost (up to a factor of 4) for both variants of the proposed controller and a broad range of tolerance. For stringent tolerance, the traditional controller marginally outperforms the proposed controller. This can be attributed to the fact that for stringent tolerances, any further reduction in step size, as prescribed by the traditional controller, leads to an increased number of time steps.

%%% ----------------------------------------------------------- %%%

\begin{figure*}
    \centering
	\includegraphics[width = 0.975\textwidth]{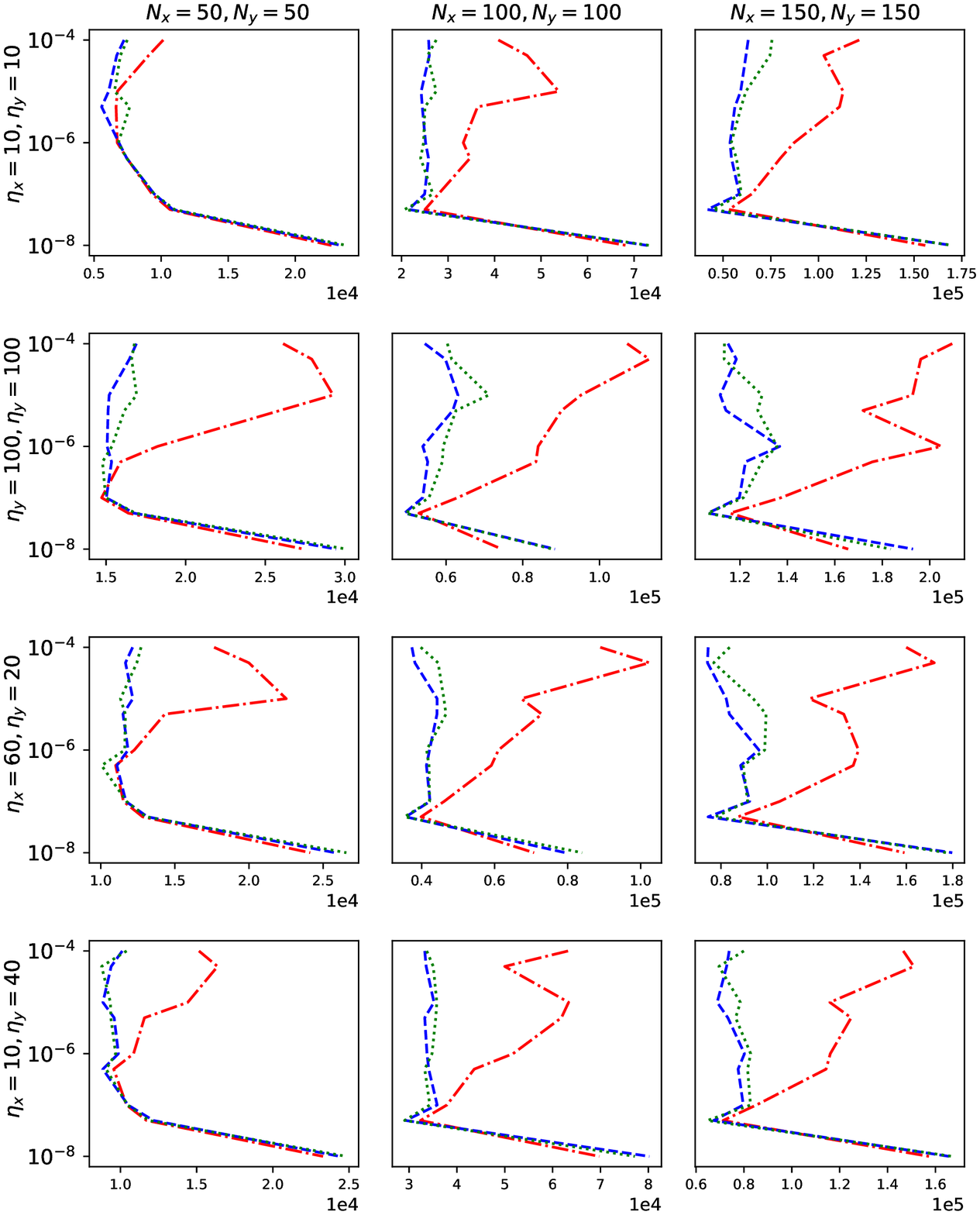}
    \caption{Figure shows the number of matrix-vector products (i.e. computational cost) vs. tolerance for the 2D porous medium equation. The values of $N_x$, $N_y$, $\eta_x$, and $\eta_y$ are varied. The red dashed-dotted lines correspond to the traditional controller, green dotted lines the penalized variant, and the blue dashed lines represent the non-penalized variant.}
    \label{fig:exprb43_porous_2d}
\end{figure*}

We, now, consider the porous medium equation in two dimensions 
\begin{equation}
    \frac{\partial u}{\partial t} = \eta_x \frac{\partial u}{\partial x} + \eta_y \frac{\partial u}{\partial y}
    + \frac{\partial^2 u^m}{\partial x^2} + \frac{\partial^2 u^m}{\partial y^2}, \nonumber
\end{equation}
where $m = 2$ and $\eta_x$ and $\eta_y$ are the components of the P\'eclet number along the $X$ and $Y$ directions, respectively. The initial condition is
\begin{equation}
    u(x, y, t = 0) = 1 + \Theta(x_1 - x) + \Theta(x - x_2) + \Theta(y_1 - y) + \Theta(y - y_2), \nonumber
\end{equation}
where $x_1 = y_1 = 0.25$, $x_2 = y_2 = 0.6$, and $\Theta$ is the Heaviside function. This initial state corresponds to a cuboid. Similar to the 1D case, the nonlinear diffusion rapidly smears out the discontinuity resulting in a progressively smoother solution. The work-precision diagram for the 2D scenario (Fig. \ref{fig:exprb43_porous_2d}) shows that the proposed step-size controller improves the performance in almost all the considered configurations. In addition, the curves are flattened-out to a large extent.

%%% ----------------------------------------------------------- %%%
%%% ----------------------------------------------------------- %%%

\subsection{Advection--Diffusion--Reaction Equation}

The final example considered is the 1D advection--diffusion--reaction (ADR) equation \cite{Tokman13}. This is similar to the combustion problem presented in \cite{Adjerid86}. The equation reads
\begin{equation*}
    \frac{\partial u}{\partial t} = \eta \frac{\partial u}{\partial x} + \frac{\partial^2 u}{\partial x^2} + \alpha u\left(u - \tfrac{1}{2}\right)(1 - u), 
\end{equation*}
where $\eta$ is the P\'eclet number and $\alpha = 1$. The initial condition is chosen to be
\begin{equation*}
	u(x, t = 0) =  256(x - x^2)^2 + 0.3.
\end{equation*}
The simulations are carried out till $t = 5\times 10^{-2}$. Different values of $\eta$ indicate varying amounts of advection with respect to diffusion, whereas the reaction rate, $\alpha$, is chosen to be a constant. The performance of the step size controllers for different number of grid points ($N$) and P\'eclet number is illustrated in Fig. \ref{fig:exprb43_adr_1d}. Consistent with the previous examples, we see that the proposed controller outperforms the traditional controller, in the lenient to moderate tolerance regime, for the high-resolution simulations. The enhanced performance becomes more prominent as the value of $\eta$ gets larger.

%%%%%%%%%%%%%%%%%%%%%%%%%%%%%%%%%%%%%%%%%%%%%%%%%%%%%%%%%%%%%%%%%%%%%%%%%%%%%%%%%%%%%%%%%%%%%%%%%%%%%%%%%%%

\section{Conclusions}
\label{sec:conclude}

In this manuscript, we have considered an adaptive step-size controller for exponential Rosenbrock integrators. The fundamental principle of this proposed step-size controller is that the step size is chosen in such a way that minimizes the computational cost, constrained by the maximum allowed time step-size; the latter being set by accuracy considerations. This allows the step-size controller to adaptively decrease the step size which can drastically improve the performance. Specifically, we have used an embedded exponential Rosenbrock integrator, EXPRB43, which has a third-order error estimate. The implementation of this time integrator involves polynomial interpolation at Leja points to compute the action of the required matrix functions. A comprehensive comparison of the proposed step-size controller with the traditional controller has been presented for different values of the P\'eclet number, the number of the grid points, and the user-specified tolerance for the various equations under consideration. We summarize our results as follows:

\begin{figure*}[th]
    \centering
	\includegraphics[width = \textwidth]{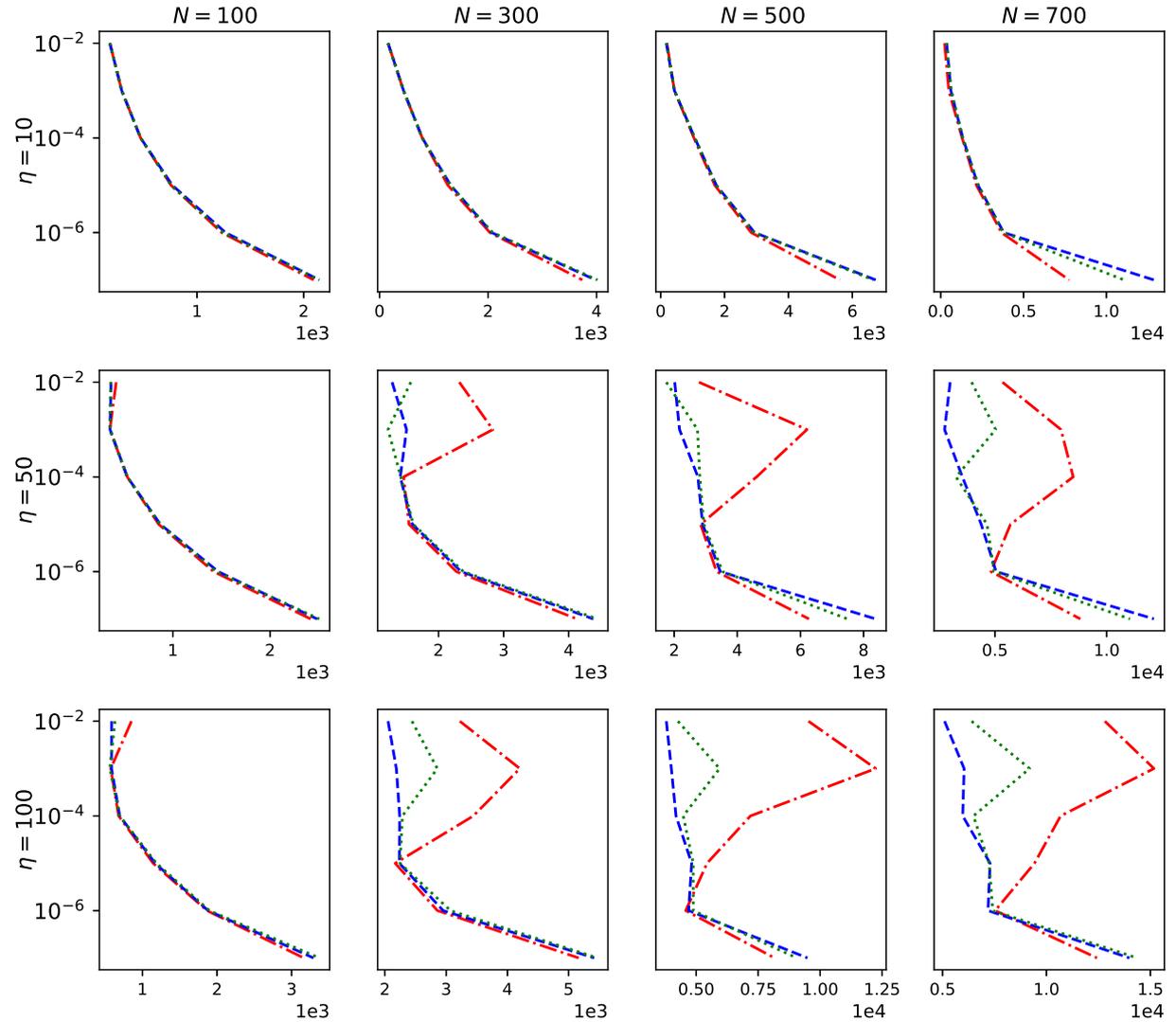}
    \caption{Figure shows the number of matrix-vector products (i.e. computational cost) vs. tolerance for the 1D ADR equation for different values of $\eta$ and $N$. The red dashed-dotted lines correspond to the traditional controller, green dotted lines the penalized variant, and the blue dashed lines represent the non-penalized variant of the proposed step-size controller.}
    \label{fig:exprb43_adr_1d}
\end{figure*}

\begin{itemize}

    \item The proposed step-size controller has superior performance (compared to the traditional controller) for almost all configurations considered here. This is particularly true in the lenient to medium tolerance regime. Arguably, equations in physics and astrophysics are solved up to accuracy within this range of tolerances. The zig-zag curves, yielded by the traditional controller, are flattened out. This exemplifies the use of such a step-size controller in practice.
    
    \item Multiple small step sizes, in many situations, do indeed incur less computational effort than a single large step size and we have seen that the proposed step-size controller can effectively exploit this fact.
    
    \item  It has been observed that the curves for the work-precision diagrams have similar `shapes' in 1D and 2D. This indicates the reliability and effectiveness of the proposed controller.
    
    \item Comparisons with explicit (RKF45) and implicit (SDIRK23) schemes have shown that the exponential Rosenbrock integrator outperforms both these classes of integrators by a significant margin.
    
    \item We recommend the non-penalized variant of the step-size controller (where the parameters have been optimized to minimize the computational cost) as it shows improved performance in almost all configurations considered. 
    
\end{itemize}

We note that one of the other highly efficient iterative schemes for exponential integrators, in the case of large and sparse matrices, is the Krylov subspace algorithm \cite{Moler03, Tokman06, Ostermann10}. The general idea of the Krylov method is to project the $\varphi_l$ function, applied to some vector, onto a Krylov subspace, of size $m$, using the Arnoldi algorithm \cite{Vorst87}. This reduces the problem of having to compute the exponential-like function of a large matrix to computing the exponential-like function of a small ($m \times m$) Hessenberg matrix. The Krylov subspace algorithm, for exponential integrators, has been shown to be highly competitive with the state-of-the-art implicit (and explicit) integrators even for highly nonlinear problems, e.g. the set of magnetohydrodynamical equations \cite{Tokman02, einkemmer2017, Deka22}. Taking advantage of the fact that any stage of an exponential integrator can be expressed in the form \[\varphi_0(A) v_0 + \varphi_1(A) v_1 + \varphi_2(A) v_2 + \hdots + \varphi_p(A) v_p, \] \cite{Niesen09} developed the \texttt{phipm} algorithm, for the Krylov method, that computes the linear combination of $\varphi_l$ functions. This is based on the idea that it is often computationally cheaper to compute the exponential of an augmented Hessenberg matrix than to evaluate several individual $\varphi_l$ functions \cite{Saad92, Sidje98, Higham11}. The work by \cite{Niesen09} also adopted a time-stepping strategy where they subdivide a given step size into several substeps, i.e. $\Delta t = \Delta t_1 + \Delta t_2 + ... + \Delta t_k$. This reduces the dimension of the Krylov subspace, and since smaller dimensional subspaces converge faster, significant computational savings can be achieved with this measure. The idea of using smaller step sizes to improve the convergence rate, and consequently the computational cost is similar to the one proposed here. The main difference being that our proposed step size controller reduces the step size of the integrator as a whole (which has the added benefit of increasing accuracy in the case of nonlinear problems), whereas the substepping in \texttt{phipm} subdivides a given step size individually for each stage of an exponential integrator. Additionally, our approach directly learns the computational cost from measurements taken during the integration as opposed to relying on, in general, an inaccurate estimate based on the sparsity structure of the matrix as is done in \texttt{phipm}. In principle, one could also combine our proposed step size controller with the \texttt{phipm} algorithm to get further savings (e.g.~in cases where different substepping strategies for the different matrix functions might be advantageous). An improvement over the \texttt{phipm} algorithm has been proposed by \cite{Gaudreault18} with the so-called Krylov with Incomplete Orthogonalization Procedure Solver (\texttt{KIOPS}) algorithm. They resort to an incomplete orthogonalization to compute the basis and present an improved Krylov adaptivity procedure. The Leja interpolation method has the advantage that it is more efficient in the parallel-computing context as only matrix-vector products are required. We are currently working on comparisons of \texttt{phipm} and \texttt{KIOPS} with the Leja method, and potentially developing a modified iterative scheme based on these methods, for a range of multi-stage exponential integrators.

This study has been performed for a set of representative but relatively simple problems. In future work, we will implement this step-size controller as part of a software package and test it for more realistic scenarios. A typical example would be the propagation of fluids or particles in the interstellar or intergalactic medium in 3D. This may include a combination of linear and/or nonlinear advection and diffusion coupled to other physical processes like dispersion, collisions, etc.

\section*{Acknowledgements}
This work is supported by the Austrian Science Fund (FWF) project id: P32143-N32. We would like to thank the two anonymous referees for their constructive criticism of this manuscript.

%%%%%%%%%%%%%%%%%%%%%%%%%%%%%%%%%%%%%%%%%%%%%%%%%%%%%%%%%%%%%%%%%%%%%%%%%%%%%%%%%%%%%%%%%%%%%%%%%%%%%%%%%%%

\appendix
\section*{Spatial Discretization}
\label{app:diff_adv}

We use the third-order upwind scheme to discretize the advective term which is given by 
\begin{equation}
    \frac{\partial u}{\partial x} \approx \frac{-u_{i + 2} + 6u_{i + 1} - 3u_{i} - 2u_{i - 1}}{6 \, \Delta x} \nonumber
\end{equation}
The primary advantage of this scheme is that it introduces less numerical diffusion. The
structure or features of the physical parameter under consideration is preserved to a large extent. The diffusive term is discretized using the second-order centred difference scheme
\begin{equation}
    \frac{\partial^2 u}{\partial x^2} \approx \frac{u_{i + 1} - 2u_i + u_{i - 1}}{\Delta x^2} \nonumber
\end{equation}

%%%%%%%%%%%%%%%%%%%%%%%%%%%%%%%%%%%%%%%%%%%%%%%%%%%%%%%%%%%%%%%%%%%%%%%%%%%%%%%%%%%%%%%%%%%%%%%%%%%%%%%%%%%

\bibliography{ref.bib}

\end{document}